\documentclass[11pt]{article}

\usepackage{latexsym}
\usepackage{amsfonts}
\usepackage{amssymb}
\usepackage{amsmath}
\usepackage{mathrsfs}
\usepackage{hyperref}

\newcommand{\be}{\begin{equation}}
\newcommand{\ee}{\end{equation}}
\newcommand{\bea}{\begin{eqnarray}}
\newcommand{\eea}{\end{eqnarray}}
\newcommand{\bean}{\begin{eqnarray*}}
\newcommand{\eean}{\end{eqnarray*}}
\newcommand{\brray}{\begin{array}}
\newcommand{\erray}{\end{array}}

\newtheorem{dfn}{Definition}[section]
\newtheorem{thm}[dfn]{Theorem}
\newtheorem{lmma}[dfn]{Lemma}
\newtheorem{ppsn}[dfn]{Proposition}
\newtheorem{crlre}[dfn]{Corollary}
\newtheorem{xmpl}[dfn]{Example}
\newtheorem{rmrk}[dfn]{Remark}

\newcommand{\bdfn}{\begin{dfn}\rm}
\newcommand{\bthm}{\begin{thm}}
\newcommand{\blmma}{\begin{lmma}}
\newcommand{\bppsn}{\begin{ppsn}}
\newcommand{\bcrlre}{\begin{crlre}}
\newcommand{\bxmpl}{\begin{xmpl}}
\newcommand{\brmrk}{\begin{rmrk}\rm}

\newcommand{\edfn}{\end{dfn}}
\newcommand{\ethm}{\end{thm}}
\newcommand{\elmma}{\end{lmma}}
\newcommand{\eppsn}{\end{ppsn}}
\newcommand{\ecrlre}{\end{crlre}}
\newcommand{\exmpl}{\end{xmpl}}
\newcommand{\ermrk}{\end{rmrk}}

\newcommand{\cla}{\mathcal{A}}
\newcommand{\clb}{\mathcal{B}}

\newcommand{\clh}{\mathcal{H}}

\newcommand{\clg}{\mathcal{G}}

\author{S. Sundar}
\title{On a construction due to Khoshkam and Skandalis }

\begin{document}
\maketitle
\begin{abstract}
In this paper, we consider the Wiener Hopf algebra, denoted $\mathcal{W}(A,P,G,\alpha)$, associated to an action of a discrete subsemigroup $P$ of a group $G$ on a $C^{*}$-algebra $A$.   We show that $\mathcal{W}(A,P,G,\alpha)$ can be represented as a groupoid crossed product. As an application, we show that when $P=\mathbb{F}_{n}^{+}$, the free semigroup on $n$ generators, the $K$-theory of $\mathcal{W}(A,P,G,\alpha)$ and the $K$-theory of $A$ coincides. 
\end{abstract}

\noindent {\bf AMS Classification No. :} {Primary 22A22; Secondary 54H20, 43A65, 46L55.}  

\section{Introduction}
\footnote{\textbf{Acknowledgements}: 
This work is done during the author's visit to the University of Orleans during May-June 2015. I  thank Prof. Renault for providing me with a visit and for discussing with me for a month at Orleans. I  thank Prof. Skandalis for spending a day at Orleans to discuss with me and Prof. Renault.  I thank him especially for explaining the idea of the proof of Theorem 6.4. I thank Prof. Philippon, IMJ Paris, for the financial support (IRSES-MODULI) which made this visit possible. 
}

Various constructions of semigroup crossed products  have been studied  by several authors most notably by Murphy and Exel. We refer to the papers \cite{Murphy94}, \cite{Murphy96} and the references therein for some of the earlier work on semigroup crossed products. We also refer the reader to \cite{Exel_endo} for Exel's work on semigroup crossed products. In view of the recent advancements in semigroup $C^{*}$-algebras (See  \cite{Li-semigroup}, \cite{Li13}), it is timely to revisit some of the constructions. 

In this paper we analyse the construction of Khoskham and Skandalis carried out in \cite{KS97} for discrete semigroups. In \cite{KS97}, actions of  the discrete semigroup $\mathbb{N}$ and of the  continuous semigroup $\mathbb{R}_{+}$ are considered.  Even though the construction due to Khoskham and Skandalis can be carried out to topological semigroups, we restrict ourselves to the discrete case. For we believe that the discrete case alone is interesting for its own sake and might be interesting to a wider audience. We plan to pursue the topological case elsewhere. We now give an overview of the problem considered and the results obtained.

 Let $P \subset G$ be a subsemigroup of a discrete group $G$ containing the identity. Let $A$ be a $C^{*}$-algebra and let $\alpha:P \to End(A)$ be a left action of $P$ on $A$ by endomorphisms. For the introduction, let us assume, for simplicity, that $A$ is unital and the action of $P$ on $A$ is by unital endomorphisms. Consider the Hilbert $A$-module $\mathcal{E}:=A \otimes \ell^{2}(P)$. For $x \in A$ and $a \in P$, define the operators $\pi(x)$ and $V_{a}$ by the equation:
\begin{align*}
\pi(x)(y \otimes \delta_{b})&:=\alpha_{b}(x)y\otimes \delta_{b}, \\
V_{a}(y \otimes \delta_{b})&:=y \otimes \delta_{ba}.
\end{align*}
Note that $V_{a}$ is an isometry for each $a \in P$ and $V_{a}V_{b}=V_{ba}$. Moreover the covariance condition $\pi(x)V_{a}=V_{a}\pi(\alpha_{a}(x))$ is satisfied for $x \in A$ and $a \in P$. It is very natural to define the "reduced crossed product" $A \rtimes_{red} P$ as the $C^{*}$-algebra generated by $\{\pi(x):x \in A\}$ and $\{V_{a}: a \in P\}$. When $A=\mathbb{C}$, then $\mathbb{C} \rtimes_{red} P$ is the reduced $C^{*}$-algebra of the semigroup $P$, denoted $C_{red}^{*}(P)$, and is studied in detail by Li in \cite{Li-semigroup} and \cite{Li13}. 

There is a related $C^{*}$-algebra called the Wiener-Hopf algebra which is much easier to understand. For $g \in G$, let $W_{g}$ be the operator on $\mathcal{E}$ defined by \begin{equation*}\label{eq:ykq}
 W_{g}(y\otimes \delta_{b})=\begin{cases}
 y \otimes \delta_{bg} & \mbox{ if
} bg \in P,\cr
   &\cr
    0 &  \mbox{ if } bg \notin P.
         \end{cases}
\end{equation*}
The $C^{*}$-algebra generated by $\{\pi(x):x \in A\}$ and $\{W_{g}: g \in G\}$ is called the Wiener-Hopf algebra associated to $(A,P,G,\alpha)$ and we denote it by $\mathcal{W}(A,P,G,\alpha)$. Clearly $A \rtimes_{red} P$ is contained in $\mathcal{W}(A,P,G,\alpha)$. Upto the author's knowledge, it is not known if $A \rtimes_{red} P=\mathcal{W}(A,P,G,\alpha)$. However the equality $A \rtimes_{red} P=\mathcal{W}(A,P,G,\alpha)$ holds if either $PP^{-1}=G$ or if the pair $(P,G)$ is quasi-lattice ordered in the sense of Nica (\cite{Nica92}). 

In this paper, we represent $\mathcal{W}(A,P,G,\alpha)$ as a groupoid crossed product. The idea of representing the Wiener-Hopf algebra as a groupoid crossed product dates back to \cite{Renault_Muhly}. In fact, we prove that $\mathcal{W}(A,P,G,\alpha)$ is isomorphic to the reduced crossed product of the form $\mathcal{D} \rtimes \clg$ where $\clg$ is the Wiener-Hopf groupoid associated to $(P,G)$. This is the content of Theorem \ref{main theorem}.  We imitate the construction of Khoskham and Skandalis carried out in \cite{KS97} for the semigroup $\mathbb{N}$ to construct the bundle $\mathcal{D}$. This forms the content of Sections 3 and 4.

In \cite{Li-Cuntz-Echterhoff}, crossed products by automorphic actions of semigroups are considered. A remarkable result that $K_{*}(A) \cong K_{*}(A \rtimes_{red} P)$ is obtained when $(P,G)$ is quasi-lattice ordered and under certain assumptions on $G$. Here we prove a very modest result that when $P=\mathbb{F}_{n}^{+}$, $K_{*}(A) \cong K_{*}(A \rtimes_{red} \mathbb{F}_{n}^{+})$. In fact the natural inclusion $A \ni x \to x \in A \rtimes_{red} \mathbb{F}_{n}^{+}$ is a $KK$-equivalence. This $KK$-equivalence relies on the amenability of the Wiener-groupoid of the pair $(\mathbb{F}_{n}^{+},\mathbb{F}_{n})$ and the fact that $\mathcal{D} \rtimes \clg$ can be expressed in terms of generators and relations. 
In Section 5, we prove that $\mathcal{D} \rtimes \clg$ can be expressed in terms of generators and relations for any quasi-lattice ordered pair $(P,G)$. The $KK$-equivalence between $A$ and $A \rtimes_{red} \mathbb{F}_{n}^{+}$ is proved in Section 6. 

Throughout this paper, we assume that the $C^{*}$-algebras bearing the name $A$ are assumed to be separable. Discrete groups are assumed to be countable.

\section{Preliminaries}
In this section, we recall the essentials regarding $C_{0}(X)$-algebras and groupoid crossed products. This is to make the paper easier to read and also to fix notations. The reader is referred to \cite{MW08} and \cite{KS04} for more on groupoid crossed products. Also the thesis \cite{Goehle} forms a good reference.

Let $G$ be a discrete group and $X$ be a right $G$-space. For $f \in C_{0}(X)$ and $g \in G$, we let $R_{g}(f) \in C_{0}(X)$ be defined by $R_{g}(f)(x)=f(xg)$ for $x \in X$.
A $C_{0}(X)$-$G$ algebra is a $C^{*}$-algebra $A$ together with a group homomorphism $\alpha:G \to Aut(A)$ and a $*$-algebra homomorphism $\rho:C_{0}(X) \to M(A)$, where $M(A)$ is the multiplier algebra of $A$, such that
\begin{enumerate}
\item[(1)] for $f \in C_{0}(X)$ and $a \in A$, $\rho(f)a=a\rho(f)$,
\item[(2)] the representation $\rho$ is non-degenerate i.e. $\overline{\rho(C_{0}(X))A}=A$, and
\item[(3)] the homomorphism $\rho$ is $G$-equivariant i.e. for $f \in C_{0}(X), g \in G$ and $a \in A$, $\alpha_{g}(\rho(f)a)=\rho(R_{g}(f))\alpha_{g}(a)$.
\end{enumerate}
Usually we omit the symbol $\rho$ and simply write $\rho(f)a$ as $fa$ for $f \in C_{0}(X)$ and $a \in A$. A $C_{0}(X)$-$G$ algebra is sometimes called an $(X,G)$-algebra. Moreover if $G$ is the trivial group, then $C_{0}(X)$-$G$ algebras are referred as $C_{0}(X)$-algebras.

Let $A$ be an $(X,G)$-algebra. For $x \in X$, define 
\begin{align*}
C_{0}(X\backslash \{x\}):&=\{f \in C_{0}(X): f(x)=0\}, \\
I_{x}:&= \overline{C_{0}(X\backslash \{x\})A}.
\end{align*}
Then $I_{x}$ is a closed two sided ideal of $A$. We denote the quotient $A/I_{x}$ by $A_{x}$. 

Similarly for a closed subset $F \subset X$, we let 
\begin{align*}
C_{0}(X\backslash F):&=\{f \in C_{0}(X): f(x)=0 \textrm{~for~}x \in F\}, ~\textrm{and} \\
I_{F}:&=\overline{C_{0}(X\backslash F)A}.
\end{align*}

Let $\displaystyle \mathcal{A}:=\coprod_{x \in X}A_{x}$. Then $\cla$ has an upper-semicontinuous bundle structure over $X$ where the topology on $\cla$ is determined by the family of sections $X \ni x \to a+I_{x} \in A_{x}$ for $a \in A$. Let $\Gamma_{0}(X,\cla)$ be the algebra of continuous sections vanishing at infinity. The map \[A \ni a \to (X \ni x \to a+I_{x} \in A_{x}) \in \Gamma_{0}(X,\cla)\]
is a $C_{0}(X)$-algebra isomorphism. We refer the reader to Appendix C of \cite{Dana-Williams} for more details, in particular to Theorems C.25 and C.26 of \cite{Dana-Williams}.

Also the transformation groupoid $X \rtimes G$ acts on the bundle $\cla$. First note that if $g \in G$, then $\alpha_{g}$ maps $I_{x.g} \to I_{x}$. Thus $\alpha_{g}$ descends to an isomorphism from $A/ I_{x.g} \to A / I_{x}$ which we denote by $\alpha_{(x,g)}$. Moreover $\alpha_{(x,g)}\alpha_{(x.g,h)}=\alpha_{(x,gh)}$. Thus $\alpha:=\{\alpha_{(x,g)}\}_{(x,g)\in X \rtimes G}$ defines an action of the groupoid $X \rtimes G$ on the bundle $\cla$. 

Now we recall the definition of groupoid crossed products. Let $\clg$ be an $r$-discrete groupoid and let $p:\cla \to \mathcal{G}^{(0)}$ be an upper semicontinuous bundle. For $x \in \clg^{(0)}$, we denote the fibre $p^{-1}(x)$ by $A_{x}$. Let $\alpha:=(\alpha_{\gamma})_{\gamma \in \clg}$ be an action of $\clg$ on $\cla$. Denote the $C^{*}$-algebra of continous sections of $\cla$ vanishing at infinity by $\Gamma_{0}(\clg^{(0)},\cla)$. Let 
\[
\Gamma_{c}(\clg,r^{*}\cla):=\{f:\clg \to \cla: \textrm{~$f$ is continuous, compactly supported and~}f(\gamma) \in A_{r(\gamma)} \textrm{~for~}\gamma \in \clg \}.
\]
The vector space $\Gamma_{c}(\clg,r^{*}\cla)$ has a $*$-algebra structure where the multiplication and the involution are given by
\begin{align*}
f*g(\gamma)&:=\sum_{r(\gamma_{1})=r(\gamma)}f(\gamma_{1})\alpha_{\gamma_{1}}(g(\gamma_{1}^{-1}\gamma)) \\
f^{*}(\gamma)&=\alpha_{\gamma}(f(\gamma^{-1}))^{*}
\end{align*}
for $f,g \in \Gamma_{c}(\clg,r^{*}\cla)$. 

For $f \in \Gamma_{c}(\clg^{(0)},\cla)$, let $\widehat{f}:\clg \to \cla$ be defined by $\widehat{f}(\gamma)=f(\gamma)$ if $\gamma \in \clg^{(0)}$ and $\widehat{f}(\gamma)=0_{r(\gamma)}$ if $\gamma \notin \clg^{(0)}$. Since we have assumed that $\clg$ is $r$-discrete, the unit space $\clg^{(0)}$ is a clopen subset of $\clg$. Thus $\widehat{f} \in \Gamma_{c}(\clg,r^{*}\cla)$. Also the inclusion $\Gamma_{c}(\clg^{(0)},\cla) \to \Gamma_{c}(\clg,r^{*}\cla)$ is a $*$-algebra homomorphism. We will consider $\Gamma_{c}(\clg^{(0)},\cla)$ as a $*$-subalgebra of $\Gamma_{c}(\clg,r^{*}\cla)$. 

A $*$-representation $\pi:\Gamma_{c}(\clg,r^{*}\cla) \to B(\clh)$ on a Hilbert space $\clh$ is said to be bounded if $\displaystyle ||\pi(f)|| \leq ||f||_{\infty}:=\sup_{x \in \clg^{(0)}}||f(x)||$ for every $f \in \Gamma_{c}(\clg^{0},\cla)$. 

For $f \in \Gamma_{c}(\clg,r^{*}\cla)$, let \[||f||_{u}:=\sup \{||\pi(f)||: \pi \textrm{~is a bounded representation of $\Gamma_{c}(\clg,r^{*}\cla)$}\}.\]
It is well known that $||.||_{u}$ is a $C^{*}$-norm on $\Gamma_{c}(\clg,r^{*}\cla)$ and the completion of $\Gamma_{c}(\clg,r^{*}\cla)$ is called the full crossed product and is denoted as $\cla \rtimes \clg$.

Next we recall the definition of the reduced crossed product. Let $x \in \clg^{(0)}$ and $\pi_{x}:A_{x} \to \mathcal{L}_{B_{x}}(\mathcal{E}_{x})$ be a representation of the fibre $A_{x}$ on the Hilbert $B_{x}$-module $\mathcal{E}_{x}$. Consider the Hilbert module $B_{x}$-module $L^{2}(\clg^{(x)},\mathcal{E}_{x})$ where $\clg^{(x)}=r^{-1}(x)$. For $f \in \Gamma_{c}(\clg,r^{*}\cla)$, let $\Lambda_{x,\pi_{x}}(f)$ be the operator on $L^{2}(\clg^{(x)},\mathcal{E}_{x})$ defined by the formula
\[
(\Lambda_{x,\pi_{x}}(f)\xi)(\gamma)=\sum_{r(\gamma_{1})=r(\gamma)=x}\pi_{x}\big(\alpha_{\gamma}(f(\gamma^{-1}\gamma_{1}))\big)\xi(\gamma_{1})
\]
for $\xi \in L^{2}(\clg^{(x)},\mathcal{E}_{x})$ and $\gamma \in \clg^{(x)}$. Then $\Lambda_{x,\pi_{x}}$ is a bounded representation of $\Gamma_{c}(\clg,r^{*}\cla)$ on the Hilbert module $L^{2}(\clg^{(x)},\mathcal{E}_{x})$. We call the representation $\Lambda_{x,\pi_{x}}$ the representation induced by the representation $\pi_{x}$.

Let $I$ be any non-empty set and $\tau:I \to \clg^{(0)}$ be a map. For $i \in I$,  let $\pi_{i}:A_{\tau(i)} \to \mathcal{L}_{B_i}(\mathcal{E}_{i})$ be a representation. We say the collection $\{\pi_{i}: i \in I\}$ is faithful if the map \[\displaystyle \Gamma_{0}(\clg^{(0)},\cla) \ni f \to \displaystyle \bigoplus_{i \in I} \pi_{i}(f(\tau(i))) \in \bigoplus_{i \in I} \mathcal{L}_{B_{i}}(\mathcal{E}_{i})\] is faithful.  

Let $\{\pi_{i}:i \in I\}$ be a faithful family. For $f \in \Gamma_{c}(\clg,r^{*}\cla)$, let \[\displaystyle ||f||_{red}:=\sup_{i \in I}||\Lambda_{\tau(i),\pi_{i}}(f)||.\]
Then $||~||_{red}$ is a norm on $\Gamma_{c}(\clg,r^{*}\cla)$ and is independent of the map $\tau:I \to \clg^{(0)}$ and also of the choice of the representations $\{\pi_{i}:i  \in I\}$. The completion of $\Gamma_{c}(\clg,r^{*}\cla)$ with respect to the reduced norm $||~||_{red}$ is called the reduced crossed product and is denoted $\cla \rtimes_{red} \clg$.

\section{Wiener-Hopf algebras associated to endomorphisms of semigroups}

Throughout this section,  let $G$  denote a   discrete countable group and let $P \subset G$  be a subsemigroup containing the identity element $e$. For a $C^{*}$-algebra $A$,  let $End(A)$  be the set of $*$-algebra homomorphisms on $A$. A map $\alpha:P \to End(A)$ is called a left action of $P$ on $A$ if  
\begin{enumerate}
\item[(1)] for $a,b \in P$, $\alpha_{a}\alpha_{b}=\alpha_{ab}$, and
\item[(2)] for $x \in A$, $\alpha_{e}(x)=x$. 
\end{enumerate}
If $A$ is unital and if $\alpha_{a}$ is unital for every $a \in P$, then we call the action $\alpha$ unital.  

\textbf{Notation:} Consider the Hilbert space $\ell^{2}(P)$. For $g \in G$, let $w_{g}$ be the partial isometry on $\ell^{2}(P)$ defined as
 \begin{equation}\label{wg}
 w_{g}(\delta_{b})=\begin{cases}
  \delta_{bg} & \mbox{ if
} bg \in P,\cr
   &\cr
    0 &  \mbox{ if } bg \notin P.
         \end{cases}
\end{equation}
Here $\{\delta_{a}:a \in P\}$ denotes the standard orthonormal basis of $\ell^{2}(P)$. 

Let $A$ be a $C^{*}$-algebra and let $\alpha:P \to End(A)$ be a left action. Consider the Hilbert $A$-module $\mathcal{E}:=A \otimes \ell^{2}(P)$. For $g \in G$, let $W_{g}$ be the operator on $\mathcal{E}$ defined by $W_{g}:=1\otimes w_{g}$. For $x \in A$, let $ \pi(x)$ be the operator on $\mathcal{E}$ defined by $\pi(x)(y \otimes \delta_{b})=\alpha_{b}(x)y \otimes \delta_{b}$

Observe the following.
\begin{enumerate}
\item[(1)] For $g \in G$, $W_{g}$ is a partial isometry. Denote the final projection $W_{g}W_{g}^{*}$ by $E_{g}$. Then $\{E_{g}: g \in G\}$ is a commuting family of projection.
\item[(2)] Observe that $W_{g}^{*}=W_{g^{-1}}$ for $g \in G$. Also $W_{g}W_{h}=W_{hg}E_{h^{-1}}$
for $g,h \in G$. By taking adjoints, we also have $W_{g}W_{h}=E_{g}W_{hg}$ for $g,h \in G$.
\item[(3)] For $a \in P$, let $V_{a}:=W_{a}$. Then $V_{a}$ is an isometry for $a \in P$. Moreover the map $P \ni a \to \mathcal{L}_{A}(\mathcal{E})$ is an anti-homomorphism i.e. $V_{a}V_{b}=V_{ba}$ and $V_{e}=Id$.  
\item[(4)] Observe that $\pi(x)V_{a}=V_{a}\pi(\alpha_{a}(x))$ for $a \in P$ and $x \in A$. Also note that $\pi(x)$ commutes with $E_{g}$ for every $x \in A$ and $g \in G$.
\end{enumerate}

We let the $C^{*}$-algebra generated by $\{\pi(x)W_{g}:x \in A, g \in G\}$ be denoted by $\mathcal{W}(A,P,G,\alpha)$. We also call the $C^{*}$-algebra $\mathcal{W}(A,P,G,\alpha)$ as the Wiener-Hopf algebra associated to the quadruple $(A,P,G,\alpha)$.  

\begin{rmrk}
\label{alternate defn}
Before proceeding further, let us note that $\mathcal{W}(A,P,G,\alpha)$ can alternatively be defined as follows. Let $\rho:A \to B(\clh)$ be a faithful representation of $A$ on the Hilbert space $\clh$. Consider the Hilbert space $\widetilde{\clh}:=\clh \otimes \ell^{2}(P)$. For $g \in G$, let $\widetilde{W}_{g}=1 \otimes w_{g}$. For $x \in \clh$, let $\widetilde{\rho}(x) \in B(\widetilde{\clh})$ be defined by $\widetilde{\rho}(x)(\xi \otimes \delta_{a})=\rho(\alpha_{a}(x))\xi \otimes \delta_{a}$. Then $\mathcal{W}(A,P,G,\alpha)$ is isomorphic to the $C^{*}$-algebra generated by $\{\widetilde{\rho}(x)\widetilde{W}_{g}:x \in A, g \in G\}$. We omit the proof of this fact. 

We only indicate that this follows by Rieffel's induction. For  we obtain a faithful representation of $\mathcal{W}(A,P,G,\alpha)$ on the Hilbert space $(A \otimes \ell^{2}(P)) \otimes_{A} \clh$ and the map $(A \otimes \ell^{2}(P))\otimes_{A} \clh \ni (x \otimes \delta_{a}) \otimes \xi \to \rho(x)\xi \otimes \delta_{a} \in \widetilde{\clh}$ is an isometric embedding. The remaining details are left to the reader. 
\end{rmrk}

When $A=\mathbb{C}$, the $C^{*}$-algebra $\mathcal{W}(A,P,G,\alpha)$ is  the usual Wiener-Hopf algebra, denoted, $\mathcal{W}(P,G)$ whose  study from the groupoid perspective was initiated in \cite{Renault_Muhly}. Our aim in this article is to perform a similar analysis for $\mathcal{W}(A,P,G,\alpha)$.

From now on, we will drop the symbol $\pi$ and simply write $x$ in place of $\pi(x)$.  There is also another related $C^{*}$-algebra which we denote by $A \rtimes_{red} P$. Let $A \rtimes_{red} P$ be the $C^{*}$-algebra generated by $\{xV_{a_{1}}^{*}V_{b_{1}}V_{a_{2}}^{*}\cdots V_{a_{n}}^{*}V_{b_{n}}:x \in A, a_{i},b_{i} \in P, n \in \mathbb{N} \}$. First we prove that $A \rtimes_{red} P \subset \mathcal{W}(A,P,G,\alpha)$. This is the content of the next lemma.
 
 \begin{lmma}
 \label{closed under multiplication}
 For $T \in \mathcal{W}(A,P,G,\alpha)$ and $g \in G$, $TW_{g}, W_{g}T \in \mathcal{W}(A,P,G,\alpha)$. Consequently $A \rtimes_{red}P \subset \mathcal{W}(A,P,G,\alpha)$. 
  \end{lmma}
 \textit{Proof.}  Since $\mathcal{W}(A,P,G,\alpha)$ is generated by $\{xW_{h}: x \geq 0, h \in G\}$ and $\{W_{g}:g \in G\}$ is $*$-closed, it  suffices to prove that for $g,h \in G$ and $x\in A$ positive, $xW_{h}W_{g}, W_{g}xW_{h} \in \mathcal{W}(A,P,G,\alpha)$. 
 
Let $x \in A$ be positive and $g, h \in G$ be given. Write $x=y^{3}$ with $y$ positive. Now note that
\begin{align*}
(yW_{h})(yW_{h})^{*}yW_{gh}&=yE_{h}y^{2}W_{gh} \\
                           & = y^{3}E_{h}W_{gh} \textrm{~~(since $E_{h}$ commutes with $y$)} \\
                           & = x W_{h}W_{g} \textrm{~~(since $W_{h}W_{g}=E_{h}W_{hg}$)}. 
\end{align*} 
Hence $xW_{h}W_{g} \in \mathcal{W}(A,P,G,\alpha)$. Note that $W_{g}xW_{h}=(yW_{g^{-1}})^{*}(y^{2}W_{h}) \in \mathcal{W}(A,P,G,\alpha)$. Thus $\mathcal{W}(A,P,G,\alpha)$ is closed under right and left multiplication by $\{W_{g}:g \in G\}$. Now the inclusion $A \rtimes_{red} P \subset \mathcal{W}(A,P,G,\alpha)$ is immediate. \hfill $\Box$

 When $A=\mathbb{C}$, $\mathbb{C} \rtimes_{red} P$  is the reduced $C^{*}$-algebra of the semigroup $P$, denoted $C_{red}^{*}(P)$ (See \cite{Li-semigroup}, \cite{Li13}). Upto the author's knowledge, even for the case $A=\mathbb{C}$, it is not known whether the equality $C_{red}^{*}(P)=\mathcal{W}(P,G)$ holds or not. We discuss two situations where the equality $A \rtimes_{red} P=\mathcal{W}(A,P,G,\alpha)$ holds.

\textbf{Nica's quasi-lattice ordered semigroups:} Recall from \cite{Nica92}, the pair $(P,G)$ is called quasi-lattice ordered if $P \cap P^{-1}=\{e\}$ and the following holds: Let $g \in G$. If $Pg \cap P$ is non-empty, then there exists $a \in P$ such that $Pg \cap P=Pa$. Let $(P,G)$ be quasi-lattice ordered and let $\alpha:P \to End(A)$ be an action of $P$ on a $C^{*}$-algebra $A$. First observe that for $g \in G$, $W_{g} \neq 0$ if and only if $Pg \cap P \neq \emptyset$. If $Pg \cap P = Pa$ then  $W_{g}=V_{a}V_{b}^{*}$ where $b=ag^{-1}$. 

Thus $\mathcal{W}(A,P,G,\alpha)$ is the $C^{*}$-algebra generated by $\{xV_{a}V_{b}^{*}:x \in A, a,b \in P\} \subset A \rtimes_{red} P$.  This proves that $\mathcal{W}(A,P,G,\alpha)=A \rtimes_{red} P$. The linear span of $\{V_{a}xV_{b}^{*}:x \in A, a,b \in P\}$ is a dense $*$-subalgebra of $\mathcal{W}(A,P,G,\alpha)$. To see this, note that if $Pa \cap Pb=\emptyset$, then $V_{b}^{*}V_{a}=0$. If $Pa \cap Pb=Pc$ then $V_{b}^{*}V_{a}=V_{t}V_{s}^{*}$ where $s,t \in P$ are such that $sa=tb=c$. Let $a_{1},b_{1},a_{2},b_{2} \in P$ and $x,y \in A$. Then $(V_{a_1}xV_{b_1}^{*})(V_{a_2}yV_{b_2}^{*})=0$ if $Pb_{1} \cap Pa_{2}=\emptyset$. If $Pb_{1} \cap Pa_{2}=Pc$ then $(V_{a_1}xV_{b_1}^{*})(V_{a_2}yV_{b_2}^{*})=V_{sa_1}\alpha_{s}(x)\alpha_{t}(y)V_{tb_2}^{*}$ where $s,t \in P$ are such that $sb_{1}=ta_{2}=c$. 
 We also note that when $P=\mathbb{N}$ and $G=\mathbb{Z}$,  $\mathcal{W}(A,\mathbb{N},\mathbb{Z},\alpha)$ agrees with the Toeplitz algebra considered in \cite{KS97}.

\begin{rmrk}
Nica in \cite{Nica92} considers left regular representations whereas we consider the right regular one. We consider here  a "right" variant of Nica's definition of a quasi lattice ordered pair.
\end{rmrk}

\textbf{Ore semigroups:}  A subsemigroup $P \subset G$ is called a right Ore subsemigroup of $G$ if $PP^{-1}=G$.
Let $P$ be a right Ore subsemigroup of $G$ and let $\alpha:P \to End(A)$ be an action of $P$ on a $C^{*}$-algebra $A$. First note that if $g=ab^{-1}$ with $a,b \in P$ then $W_{g}=V_{b}^{*}V_{a}$.
Hence $\mathcal{W}(A,P,G,\alpha)$ is the $C^{*}$-algebra generated by $\{xV_{a}^{*}V_{b}: a, b \in P\} \subset A \rtimes_{red} P$. This implies that $\mathcal{W}(A,P,G,\alpha)=A \rtimes_{red} P$.

We consider now a unitisation procedure. Let $A$ be a $C^{*}$-algebra and $\alpha:P \to End(A)$ be a left action. Let $A^{+}:=A \oplus \mathbb{C}$ be the unitisation of $A$. For $a \in P$, let $\alpha^{+}_{a}:A^{+} \to A^{+}$ be defined by $\alpha^{+}_{a}(x,\lambda)=(\alpha_{a}(x),\lambda)$. Then $\alpha^{+}: P \to End(A^{+})$ is a unital left action. We call $(A^{+},P,G,\alpha^{+})$ as the unitisation of $(A,P,G,\alpha)$. 

\begin{lmma}
\label{unitisation}
Let $A$ be a $C^{*}$-algebra and $\alpha: P \to End(A)$ be a left action. We have the following short exact sequence
\[
0 \rightarrow \mathcal{W}(A,P,G,\alpha) \rightarrow \mathcal{W}(A^{+},P,G,\alpha^{+}) \rightarrow \mathcal{W}(P,G) \rightarrow 0.
\]
\end{lmma}
\textit{Proof.} Let $\rho$ be a faithful representation of $A^{+}$ on a Hilbert space $\clh$. Then by Remark \ref{alternate defn}, it follows that $\mathcal{W}(A^{+},P,G,\alpha^{+})$ is the $C^{*}$-algebra generated by $\{\widetilde{\rho}(x):x \in A\}$ and $\{\widetilde{W}_{g}: g \in G$ where for $x \in A$, $\widetilde{\rho}(x)$ and $\widetilde{W}_{g}$ are defined as in Remark \ref{alternate defn}. Since $\rho$ is faithful on $A$, it follows by Remark \ref{alternate defn} that $\mathcal{W}(A,P,G,\alpha)$ is the $C^{*}$-algebra generated by $\{\widetilde{\rho}(x)\widetilde{W}_{g}:x \in A, g \in G\}$. Clearly $\mathcal{W}(A,P,G,\alpha)$ is an ideal in $\mathcal{W}(A^{+},P,G,\alpha^{+})$.

Now consider the Hilbert $A^{+}$-module $\mathcal{E}^{+}:=A^{+} \otimes \ell^{2}(P)$. For $g \in G$, let $W_{g}=1\otimes w_{g}$. For $x \in A$, let $\pi^{+}(x)$ be the operator on $\mathcal{E}^{+}$ defined by \[\pi^{+}(x)((y,\lambda)\otimes \delta_{a})=(\alpha_{a}^{+}(x)y+\lambda \alpha_{a}^{+}(x),0) \otimes \delta_{a}.\]
Then the argument in the preceeding paragraph implies that $\mathcal{W}(A,P,G,\alpha)$ is the $C^{*}$-algebra generated by $\{\pi^{+}(x)W_{g}: x \in A,g \in G\}$ and is an ideal in $\mathcal{W}(A^{+},P,G,\alpha)$. By definition, $\mathcal{W}(A^{+},P,G,\alpha)$ is the $C^{*}$-algebra generated by $\{\pi^{+}(x): x \in A\}$ and $\{W_{g}:g \in G\}$. 

Denote the map $A^{+} \ni (x,\lambda) \to \lambda \in \mathbb{C}$ by $\epsilon$. Clearly $\mathcal{E}^{+} \otimes_{\epsilon} \mathbb{C} \cong \ell^{2}(P)$. Denote the map $\mathcal{L}_{A^{+}}(\mathcal{E}^{+}) \ni T \to T \otimes 1 \in B(\mathcal{E}^{+} \otimes_{\epsilon} \mathbb{C}) \cong B(\ell^{2}(P))$ by $\widetilde{\epsilon}$. Then clearly $\widetilde{\epsilon}$ vanishes on $\mathcal{W}(A,P,G,\alpha)$ and the $C^{*}$-algebra $\mathcal{W}(A^{+},P,G,\alpha^{+})$ is mapped onto the Wiener-Hopf algebra $\mathcal{W}(P,G)$. 
 Let $\sigma:\mathcal{W}(P,G) \to \mathcal{W}(A^{+},P,G,\alpha^{+})$ be the $*$-homomorphism such that $\sigma(w_{g})=W_{g}$ for every $g \in G$. Observe that $\widetilde{\epsilon} \circ \sigma = id$.  

We claim that the kernel of the map $\widetilde{\epsilon}:\mathcal{W}(A^{+},P,G,\alpha^{+}) \to \mathcal{W}(P,G)$ is $\mathcal{W}(A,P,G,\alpha)$. With our notations $\mathcal{W}(A^{+},P,G,\alpha^{+})$ is generated by $\{\pi^{+}(x):x \in A\}$ and $\sigma(\mathcal{W}(P,G))$. Note that $\widetilde{\epsilon}$ vanishes on $\pi^{+}(x): x \in A \}$ and $\widetilde{\epsilon}$ is one-one on $\sigma(\mathcal{W}(P,G))$. Thus the kernel of $\widetilde{\epsilon}:\mathcal{W}(A^{+},P,G,\alpha^{+}) \to \mathcal{W}(P,G)$ is the ideal in $\mathcal{W}(A^{+},P,G,\alpha^{+})$ generated by $\{\pi^{+}(x): x \in A\}$ which is $\mathcal{W}(A,P,G,\alpha)$. This completes the proof. \hfill $\Box$

We end this section by discussing the Wiener-Hopf groupoid considered in \cite{Renault_Muhly}. The Wiener-Hopf groupoid will play a prominent role in our analysis of the algebra $\mathcal{W}(A,P,G,\alpha)$. 

Let $\mathcal{C}(G)$ denote the set of all subsets of $G$. We identify $\mathcal{C}(G)$ with the product $\{0,1\}^{G}$ by identifying a subset of $G$ with its characteristic function. We endow $\{0,1\}^{G}$ and thus $\mathcal{C}(G)$ with the product topology. Endowed with this product topology, $\mathcal{C}(G)$ is a compact Hausdorff space. Also $G$ acts on $\mathcal{C}(G)$ by translation. The action of $G$ on $\mathcal{C}(G)$ is given by $\mathcal{C}(G) \times G \ni (A,g) \to Ag \in \mathcal{C}(G)$. We denote the closure of $\{P^{-1}a: a \in P\}$ in $\mathcal{C}(G)$ by $\Omega$ and we call it the Wiener-Hopf compactification of the pair $(P,G)$. Note that $\Omega$ is compact. Let \[\widetilde{\Omega}:=\{Ag: A \in \Omega, g \in G\}.\] Then $\widetilde{\Omega}$ is $G$-invariant.

The basic facts concerning the spaces $\widetilde{\Omega}$ and $\Omega$ are summarised in the next lemma.

\begin{lmma}
\label{facts on Omega}
With the foregoing notations, we have the following.

\begin{enumerate}
\item[(1)] Let $A \in \widetilde{\Omega}$ and $g \in G$. Then $Ag \in \Omega$ if and only if $g^{-1} \in A$.
\item[(2)] The set $\Omega$ is a compact open subset of $\widetilde{\Omega}$.
\item[(3)] The space $\widetilde{\Omega}$ is locally compact.
\item[(4)] The $C^{*}$-subalgebra generated by $\{1_{\Omega g}: g \in G\}$ is dense in $C_{0}(\widetilde{\Omega})$.
\end{enumerate}

\end{lmma}
\textit{Proof.} Observe that if $B \in \Omega$ then, it follows from the definition that, $e \in B$ and $P^{-1}B \subset B$. Now let $A \in \widetilde{\Omega}$ and $g \in G$ be given. Suppose $Ag \in \Omega$. Since $e \in Ag$, it follows that $g^{-1} \in A$. Now assume $g^{-1} \in A$. Write $A=Bh$ with $B \in \Omega$ and $h \in G$. Let $(a_{n})$ be a sequence in $P$ such that $P^{-1}a_{n} \to B$. Now $g^{-1}h^{-1} \in B$. Thus there exists $N \in \mathbb{N}$ such that $g^{-1}h^{-1}a_{n}^{-1} \in P^{-1}$ for $n \geq N$. For $n \geq N$, let $b_{n} \in P$ be such that $b_{n}=a_{n}hg$. Then $P^{-1}b_{n} \to Bhg=Ag$. Thus $Ag \in \Omega$. This proves $(1)$.

Observe that $(1)$ implies that $\Omega=\{A \in \widetilde{\Omega}: 1_{A}(e)=1\}$. Thus $\Omega$ is open in $\widetilde{\Omega}$. By definition, $\Omega$ is compact. This proves $(2)$.  The local compactness of $\widetilde{\Omega}$ follows from $(2)$ and by the equality $\widetilde{\Omega}=\bigcup_{g \in G}\Omega g$. 

Note that by $(1)$, for $A \in \widetilde{\Omega}$, $1_{\Omega g}(A)=1_{A}(g)$. The last assertion is then an immediate consequence of the Stone-Weierstrass theorem.  This completes the proof \hfill $\Box$.

Consider the transformation groupoid $\widetilde{\Omega} \rtimes G$. The Wiener-Hopf groupoid (\cite{Renault_Muhly}), let us denote by $\mathcal{G}$, is defined as the restriction groupoid $\widetilde{\Omega} \rtimes G|_{\Omega}$.  That is \begin{align*}
\mathcal{G}:&=\{(A,g) \in \widetilde{\Omega} \rtimes G: A \in \Omega, Ag \in \Omega\} \\
            &=\{(A,g): A \in \Omega, g^{-1} \in A\}
            \end{align*}
            where the groupoid multiplication and the inversion are defined as
            \begin{align*}
            (A,g)(Ag,h)&=(A,gh)\\
            (A,g)^{-1}&=(Ag,g^{-1})
            \end{align*}
The topology on $\clg$ is the subspace topology inherited from the product topology on $\widetilde{\Omega} \times G$.  Note that $\clg$ is an $r$-discrete groupoid. Observe that for an open subset $U \subset \widetilde{\Omega}$ and $g \in G$, \[r((U \times \{g\})\cap \clg)= \{A \in U \cap \Omega: 1_{A}(g^{-1})=1\}.\]
    Thus the range map $r$ is an open map. Clearly $r$ restricted to $(U \times \{g\}) \cap \clg)$ is $1$-$1$. Hence  $r$ is a local homeomorphism.

\section{A groupoid crossed product representation of $\mathcal{W}(A,P,G,\alpha)$}

Let $A$ be a $C^{*}$-algebra and $\alpha:P \to End(A)$ be a left action. Let $\Omega$ be the Wiener-Hopf compactification of $(P,G)$ and $\widetilde{\Omega}=\bigcup_{g \in G} \Omega g$. Denote the Wiener-Hopf groupoid $\widetilde{\Omega} \rtimes G|_{\Omega}$ by $\mathcal{G}$. In this section, we show that the Wiener-Hopf algebra $\mathcal{W}(A,P,G,\alpha)$ is isomorphic to a reduced crossed product of the form $\mathcal{D} \rtimes_{red} \mathcal{G}$. To that effect, we imitate the construction due to Khoshkam and Skandalis carried out in \cite{KS97} for the semigroup $\mathbb{N}$.

We let $A^{+}:=A\oplus \mathbb{C}$ be the unitisation of $A$. Let $\ell^{\infty}(G,A^{+})$ be the $C^{*}$-algebra of bounded functions on $G$ taking values in  $A^{+}$.  The group $G$ acts on $\ell^{\infty}(G,A^{+})$ by translation. Let us denote the action of $G$ on $\ell^{\infty}(G,A^{+})$ by $\beta$. The action $\beta$ is given by the formula: for $g \in G$, $f \in \ell^{\infty}(G,A^{+})$ and $s \in G$, $\beta_{g}(f)(s)=f(sg)$.

For $g \in G$, let $j_{g}:A \to \ell^{\infty}(G,A^{+})$ be the $*$-homomorphism defined as follows: for $x \in A$ and $h \in G$, \begin{equation*}
 j_{g}(x)(h):=\begin{cases}
 \alpha_{hg^{-1}}(x) & \mbox{ if
} h \in Pg,\cr
   &\cr
    0 &  \mbox{ if } h \notin Pg.
         \end{cases}
\end{equation*}
Note that for $s \in G$, $g \in G$ and $x \in A$, $\beta_{s}(j_{g}(x))=j_{gs^{-1}}(x)$.

For $\phi \in C_{0}(\widetilde{\Omega})$, let $\widehat{\phi} \in \ell^{\infty}(G,A^{+})$ be defined by $\widehat{\phi}(g)=\phi(P^{-1}g)$ for $g \in G$. Note that the map $C_{0}(\widetilde{\Omega}) \ni \phi \to \widehat{\phi} \in \ell^{\infty}(G,A^{+})$ is $G$-equivariant. Since $\{P^{-1}a: a \in P\}$ is dense in $\Omega$, it follows that $\{P^{-1}g: g \in G\}$ is dense in $\widetilde{\Omega}$. Thus  the $*$-algebra homomorphism $ C_{0}(\widetilde{\Omega}) \ni \phi \to \widehat{\phi} \in \ell^{\infty}(G,A^{+})$ is $1$-$1$. Henceforth we will identify $C_{0}(\widetilde{\Omega})$ as a $*$-subalgebra of $\ell^{\infty}(G,A^{+})$ and will simply denote $\widehat{\phi}$ as $\phi$ for $\phi \in C_{0}(\widetilde{\Omega})$.

Let $\widetilde{D}_{0}$ be the $C^{*}$-subalgebra of $\ell^{\infty}(G,A^{+})$ generated by  $\{j_{g}(x): x \in A, g \in G\}$ and 
$C_{0}(\widetilde{\Omega})$. Denote the $C^{*}$-algebra generated by $\{j_{g}(x):x \in A, g \in G\}$ by $\widetilde{D}_{1}$ and the $C^{*}$-algebra generated by $\{\phi j_{g}(x): \phi \in C_{0}(\widetilde{\Omega}), g \in G, x \in A\}$ by $\widetilde{D}$. Observe that $\widetilde{D} \subset \ell^{\infty}(G,A)$. Moroever the commutative algebra $C_{0}(\widetilde{\Omega})$ is contained in the center of $\widetilde{D}_{0}$. Also $\widetilde{D}$ is an ideal in $\widetilde{D}_{0}$. Thus $\widetilde{D}$ is a $C_{0}(\widetilde{\Omega})$-algebra. 

The action of $C_{0}(\widetilde{\Omega})$ is given by left multiplication. To see that the action is non-degenerate, note that by $(1)$ of Lemma \ref{facts on Omega} $1_{\Omega g}j_{g}(x)=j_{g}(x)$ for $g \in G$ and $x \in A$. Hence $1_{\Omega g}\phi j_{g}(x)=\phi j_{g}(x)$ for $\phi \in C_{0}(\widetilde{\Omega})$, $g \in G$ and $x \in A$. Hence the action of $C_{0}(\widetilde{\Omega})$ on $\widetilde{D}$ is non-degenerate. Moreover the equality $1_{\Omega g}j_{g}(x)=j_{g}(x)$ implies that $\widetilde{D}_{1}$ is contained in $\widetilde{D}$.

Note that $\widetilde{D}$, $\widetilde{D}_{0}$ and $\widetilde{D}_{1}$ are $G$-invariant under the translation action $\beta$. Thus $\widetilde{D}$ is a $(\widetilde{\Omega},G)$-algebra. Denote the corresponding upper-semicontinuous bundle  on which the transformation groupoid $\widetilde{\Omega} \rtimes G$ acts by $\widetilde{\mathcal{D}}$. Denote the action of $\widetilde{\Omega} \rtimes G$ on $\widetilde{\mathcal{D}}$ by $\widetilde{\beta}:=(\beta_{(X,g)})$. Recall that the fibre $\widetilde{D}_{X}$ at a point $X \in \widetilde{\Omega}$ is given by $\widetilde{D}/I_{X}$ where $I_{X}=C_{0}(\widetilde{\Omega}\backslash \{X\})\widetilde{D}$. 

Let $\mathcal{D}:=\displaystyle \coprod_{X \in \Omega} \widetilde{D}/I_{X}$ be the restriction of $\widetilde{\mathcal{D}}$ onto the subset $\Omega$.    Then the action $\widetilde{\beta}$ restricts to an action $\beta:=(\beta_{(X,g)})$ of the Wiener-Hopf groupoid $\widetilde{\Omega} \rtimes G|_{\Omega}$. Recall that the map $\beta_{(X,g)}: \widetilde{D}/I_{X.g} \to \widetilde{D}/I_{X}$ is given by $\beta_{(X,g)}(d+I_{X.g})=\beta_{g}(d)+I_X$ for $(X,g) \in \mathcal{G}$. We use the same letter $\beta$ to denote  the action of $G$ on $\ell^{\infty}(G,A^{+})$ and the action of the Wiener-Hopf groupoid on $\mathcal{D}$. We call the groupoid dynamical system $(\mathcal{D},\clg,\beta)$ the \textbf{Wiener-Hopf groupoid dynamical system} associated to the quadruple $(A,P,G,\alpha)$.

Let us make a few observations when the algebra $A$ is unital and the action $\alpha:P \to End(A)$ is unital.

\begin{rmrk}
\label{unital case}
Let $A$ be a unital $C^{*}$-algebra with the multiplicative unit $1_{A}$ and $\alpha:P \to End(A)$ be a unital action.  
Denote the map $A^{+} \ni (a,\lambda) \to a+\lambda 1_{A} \in A$ by $\delta$. Then $\delta$ is a $*$-homomorphism. Let $\widetilde{\delta}:\ell^{\infty}(G,A^{+}) \to \ell^{\infty}(G,A)$ be definded by $\widetilde{\delta}(\phi)=\delta \circ \phi$.  
\begin{enumerate}
\item The map $\widetilde{\delta}$ restricted to $\widetilde{D}$ is injective. For $\widetilde{\delta}$ is 1-1 on the ideal $\ell^{\infty}(G,A)$ contained in $\ell^{\infty}(G,A^{+})$ and $\widetilde{D} \subset \ell^{\infty}(G,A)$.
\item For $\phi \in C_{0}(\widetilde{\Omega})$, let $\overline{\phi} \in \ell^{\infty}(G,A)$ be defined by $\overline{\phi}(g)=\phi(P^{-1}g)1_{A}$ for $g \in G$. In otherwords,  $\overline{\phi}=\widetilde{\delta}(\widehat{\phi})$. Note that the map $C_{0}(\widetilde{\Omega}) \ni \phi \to \overline{\phi} \in \ell^{\infty}(G,A)$ is 1-1. 
\item Note that, by Lemma \ref{facts on Omega}, for $g \in G$, $\overline{1_{\Omega g}}=\widetilde{\delta}(j_{g}(1_{A}))$. Thus by Part (4) of Lemma \ref{facts on Omega}, it follows that $\widetilde{\delta}(\widetilde{D})$ is generated by $\{\widetilde{\delta}(j_{g}(x)):g \in G,x \in A\}$ i.e. $\widetilde{\delta}(\widetilde{D})=\widetilde{\delta}(\widetilde{D}_{1})$.
\end{enumerate}
Thus in the unital case, we will suppress the notation $\widetilde{\delta}$ and simply denote $\widetilde{\delta}(j_{g}(x))$ by $j_{g}(x)$ and $\widetilde{\delta}(\widetilde{D})$ by $\widetilde{D}$. Moreover we will denote $j_{g}(1_{A})=\overline{1_{\Omega g}}$ by $1_{\Omega g}$ and consider $C_{0}(\widetilde{\Omega})$ as a subalgebra of $\widetilde{D}$. 

In short, in the unital case, it is not necessary to pass to the unitisation $A^{+}$.

\end{rmrk}

\begin{rmrk}
\label{sections}
Note that the map \[\widetilde{D} \ni d \to (\Omega \ni X \to d+I_{X} \in \widetilde{D}/I_{X}) \in \Gamma(\Omega, \mathcal{D})\] descends to an isomorphism between $\widetilde{D}/I_{\Omega}$ and $\Gamma(\Omega,\mathcal{D})$ where $I_{\Omega}$ is the closure of $C_{0}(\widetilde{\Omega}\backslash \Omega)\widetilde{D}$. Since $\Omega$ is compact, we will simply denote $\Gamma_{0}(\Omega, \mathcal{D})$ by $\Gamma(\Omega,\mathcal{D})$. 
\end{rmrk}

Now we can state our main theorem.

\begin{thm}
\label{main theorem}
With the foregoing notations, the $C^{*}$-algebra $\mathcal{W}(A,P,G,\alpha)$ is isomorphic to the reduced groupoid crossed product $\mathcal{D} \rtimes_{red} \clg$. 
\end{thm}

First we fix some notations. For $g \in G$ and $d \in \widetilde{D}$, let $W_{d,g} \in \Gamma_{c}(\clg,r^{*}\mathcal{D})$ be defined by
\begin{equation*}
 W_{d,g}(X,h):=\begin{cases}
 d+I_{X} & \mbox{ if
} h =g,\cr
   &\cr
    0+I_{X} &  \mbox{if } h \neq g
         \end{cases}
\end{equation*}
for $(X,h) \in \clg$. Note that the linear span of $\{W_{d,g}:d \in \widetilde{D}, g \in G\}$ is $\Gamma_{c}(\clg,r^{*}\mathcal{D})$. Observe  that for $g \in G$, the map $\widetilde{D} \ni d \to W_{d,g} \in \Gamma_{c}(\clg,r^{*}\mathcal{D})$ is linear. With this notation, the subalgebra $\Gamma(\Omega,\mathcal{D}) \subset \Gamma_{c}(\clg,r^{*}\mathcal{D})$ is $ \{W_{d,e}:d \in \widetilde{D}\}$.

For $a \in P$, let $\widetilde{\pi}_{a}:\widetilde{D} \to A$ be the map defined by $\widetilde{\pi}_{a}(d)=d(a)$. Clearly $\widetilde{\pi}_{a}$ vanishes on $C_{0}(\widetilde{\Omega}\backslash \{P^{-1}a\})\widetilde{D}$ and thus descends to a $*$-homomorphism from $\widetilde{D}/I_{P^{-1}a}$ to $A$. We denote the resulting map  by $\pi_a$. 

\textit{Claim:} The map $\Gamma(\Omega,\mathcal{D}) \ni f \to \displaystyle \bigoplus_{a \in P}\pi_{a}(f(P^{-1}a)) \in \displaystyle \bigoplus_{a \in P} A$ is injective. By Remark \ref{sections}, it is enough to show that the kernel of the map $\widetilde{D} \ni d \to (d(a))_{a \in P} \in \displaystyle  \bigoplus_{a \in P}A$ is $I_{\Omega}$. Let $d \in \widetilde{D}$ be such that $d(a)=0$ for every $a \in P$. Let $\epsilon>0$ be given. Choose $\phi \in C_{0}(\widetilde{\Omega})$ such that $||d-\phi d|| \leq \epsilon$.  Since $d(a)=0$ for every $a \in P$, it follows that $1_{\Omega}d=0$. Now $\phi d = \phi(1-1_{\Omega})d$. Clearly $\phi(1-1_{\Omega}) \in C_{0}(\widetilde{\Omega}\backslash \Omega)$.  Thus $\phi d \in I_{\Omega}$. Since $I_{\Omega}$ is a closed ideal, it follows that $d \in I_{\Omega}$. This proves the claim. Let us isolate the just proved fact in a remark for latter purposes.

\begin{rmrk}
\label{the kernel}
The map $\widetilde{D}/I_{\Omega} \ni d+I_{\Omega} \to (d(a))_{a \in P} \in \ell^{\infty}(P,A)$ is  injective.
\end{rmrk}

By Lemma \ref{facts on Omega}, it follows that $\mathcal{G}^{X}=r^{-1}(X)=\{(P^{-1}a,a^{-1}b):b \in P\}$ for $X=P^{-1}a$. Recall from Section 1, the representation $\Lambda_{P^{-1}a,\pi_a}$ of $\Gamma_{c}(\mathcal{G},r^{*}\mathcal{D})$ on the Hilbert $A$-module $\ell^{2}(\mathcal{G}^{X},A) \cong \ell^{2}(a^{-1}P) \otimes A$ induced by the representation $\pi_{a}$. For $a \in P$,  we will denote the representation $\Lambda_{P^{-1}a,\pi_{a}}$ simply by $\Lambda_{a}$.

A direct verification yields that for $d \in \widetilde{D}$, $g \in G$, the operator $\Lambda_{a}(W_{d,g})$ on $\ell^{2}(a^{-1}P)\otimes A$ is given by the formula:

\begin{equation}
\label{expression}
\Lambda_{a}(W_{d,g})(\delta_{a^{-1}b}\otimes y)=\begin{cases}
 \delta_{a^{-1}bg^{-1}}\otimes d(bg^{-1})y & \mbox{ if
} b \in Pg\cr
   &\cr
    0 &  \mbox{if } b \notin Pg.
         \end{cases}
\end{equation}
We leave the verification of the above expression to the reader.

For $a \in P$, let $U_{a}:\ell^{2}(P) \otimes A \to \ell^{2}(a^{-1}P) \otimes A$ be the unitary defined by the equation $U_{a}(\delta_{b} \otimes y):=\delta_{a^{-1}b}\otimes y$. For $d \in \ell^{\infty}(G,A^{+})$, let $\pi(d)$ be the 'multiplication' operator on $\ell^{2}(P) \otimes A$ defined by \[\pi(d)(\delta_{b} \otimes y):=\delta_{b} \otimes d(b)y.\]
Then $\pi$ is a representation of $\ell^{\infty}(G,A^{+})$ on the Hilbert $A$-module $\ell^{2}(P) \otimes A$. Then Eq.\ref{expression} implies that for  $a \in P$, 
\begin{equation}
\label{norm}
U_{a}^{*}\Lambda_{a}(W_{d,g})U_{a}=\pi(d)W_{g}^{*}
\end{equation} for $d \in \widetilde{D}$ and $g \in G$. Here $W_{g}=w_{g} \otimes 1$ where $w_{g}$ stands for the operator defined by  Eq.\ref{wg}. 

\textit{Proof of Theorem \ref{main theorem}:} Since the map $\Gamma(\Omega, \mathcal{D}) \ni f \to \pi_{a}(f(P^{-1}a)) \in \displaystyle \bigoplus_{a \in P} A$ is injective, the reduced norm on $\Gamma_{c}(\clg,r^{*}\mathcal{D})$ is given by $\displaystyle ||f||_{red}=\sup_{a \in P}||\Lambda_{a}(f)||$. But $\Gamma_{c}(\clg,r^{*}\cla)$ is the linear span of $\{W_{d,g}:d \in \widetilde{D},g \in G\}$. Thus if $\displaystyle f=\sum_{i}W_{d_{i},g_{i}}$ is a finite sum, then by Eq.\ref{norm}, we have that $||f||_{red}=||\displaystyle \sum_{i}\pi(d_{i})W_{g_{i}}^{*}||$. This implies that the reduced crossed product $\mathcal{D} \rtimes_{red} \clg$ is isomorphic to the $C^{*}$-subalgebra of $\mathcal{L}_{A}(\ell^{2}(P) \otimes A)$ generated by $\{\pi(d)W_{g}:d \in \widetilde{D}, g \in G\}$. 

We claim that the $C^{*}$-algebra generated by $\{\pi(d)W_{g}:d \in \widetilde{D},g \in G\}$ is $\mathcal{W}(A,P,G,\alpha)$. Let $\mathcal{B}$ stands for the $C^{*}$-algebra generated by $\{\pi(d)W_{g}:d \in \widetilde{D},g \in G\}$. Note that for $x \in A$ and $g \in G$, $xW_{g}=\pi(j_{e}(x))W_{g}$. Hence $\mathcal{W}(A,P,G,\alpha) \subset \clb$. 

Observe that for $g \in G$ and $x \in A$, $\pi(1_{\Omega g})=E_{g}:=W_{g}W_{g}^{*}$ and $\pi(j_{g}(x))=W_{g}xW_{g}^{*}$. 
Since $C_{0}(\widetilde{\Omega})$ is generated by $\{1_{\Omega g}: g \in G\}$, it follows that the linear span of  products of the form $1_{\Omega g_{1}}1_{\Omega g_{2}} \cdots 1_{\Omega g_{m}}j_{h_1}(x_1)\cdots j_{h_{n}}(x_{n})$ forms a dense $*$-subalgebra of $\widetilde{D}$. 
Let $d:=1_{\Omega g_{1}}1_{\Omega g_{2}} \cdots 1_{\Omega g_{m}}j_{h_1}(x_1)\cdots j_{h_{n}}(x_{n}) \in \widetilde{D}$ be such a product. Then \[\pi(d)=E_{g_1}E_{g_{2}}\cdots E_{g_{m}}(W_{h_{1}}x_{1}W_{h_{1}}^{*})\cdots (W_{h_{n}}x_{n}W_{h_{n}}).\] A repeated application of Lemma \ref{closed under multiplication} implies that $\pi(d) \in \mathcal{W}(A,P,G,\alpha)$. Hence the image of $\widetilde{D}$ under $\pi$ is contained in $\mathcal{W}(A,P,G,\alpha)$. Another application of Lemma \ref{closed under multiplication} implies that $\{\pi(d)W_{g}:d \in \widetilde{D}, g \in G \} \subset \mathcal{W}(A,P,G,\alpha)$. Hence $\mathcal{B} \subset \mathcal{W}(A,P,G,\alpha)$. This completes the proof. \hfill $\Box$.  

We end this section by recording the relations satsified by $\{W_{d,g}: d \in \widetilde{D}, g \in G\}$ in the following lemma for latter purposes.

\begin{lmma}
\label{Relations}
We have the following.
\begin{enumerate}
\item[(1)] For $d \in \widetilde{D}$ and $g \in G$, $W_{d,g}=0$ if and only if $1_{\Omega g^{-1}}1_{\Omega}d=0$. Thus for every $d \in \widetilde{D}$ and $g \in G$, $W_{d,g}=W_{1_{\Omega}1_{\Omega g^{-1}}d,g}$.
\item[(2)] For $d \in \widetilde{D}$ and $g \in G$, $W_{d,g}^{*}=W_{\beta_{g}^{-1}(d^{*}),g^{-1}}$.
\item[(3)] For $d_{1},d_{2} \in \widetilde{D}$ and $g_{1},g_{2} \in G$, $W_{d_1,g_1}W_{d_2,g_2}=W_{d,g_1g_2}$ where $d=1_{\Omega g_{1}^{-1}}d_{1} \beta_{g_1}(d_{2})$.
\item[(4)] The map $\widetilde{D} \ni d \to W_{d,e} \in \mathcal{D} \rtimes \clg$ is a $*$-homomorphism whose kernel is $I_{\Omega}$. 
\end{enumerate}
\end{lmma} 
\textit{Proof.} Let $d \in \widetilde{D}$ and $g \in G$ be given. Observe that \[1_{\Omega g^{-1}}1_{\Omega}d+I_{X}=1_{\Omega g^{-1}}(X) 1_{\Omega}(X)(d+I_{X})\] for every $X \in \widetilde{\Omega}$. Now $W_{d,g}=0$ if and only if $d+I_{X}=0$ whenever $X \in \Omega$ and $Xg \in \Omega$. In other words, $W_{d,g}=0$ if and only if $1_{\Omega g^{-1}}(X)1_{\Omega}(X)(d+I_{X})=0$ for every $X \in \widetilde{\Omega}$. Hence $W_{d,g}=0$ if and only if $1_{\Omega g^{-1}}1_{\Omega} d=0$. This proves $(1)$. The second assertion of $(1)$ follows from the fact that for $g \in G$, the map $\widetilde{D} \ni d \to W_{d,g} \in \Gamma_{c}(\clg,r^{*}\mathcal{D})$ is linear. 

We leave the verification of $(2)$ to the reader.

Let $d_{1}, d_{2} \in \widetilde{D}$ and $g_{1},g_{2} \in G$ be given. Observe that for $(X,g) \in \mathcal{G}$, 
\begin{align*}
W_{d_{1},g_{1}}W_{d_{2},g_{2}}(X,g)&=\displaystyle \sum_{(X,s) \in \mathcal{G}}W_{d_1,g_1}(X,s)\beta_{(X,s)}(W_{d_2,g_2}(Xs,s^{-1}g)) \\
                                   &= \sum_{s \in G}1_{\Omega s^{-1}}(X)W_{d_1,g_1}(X,s)\beta_{(X,s)}(W_{d_2,g_2}(Xs,s^{-1}g)).
\end{align*}
The above equation implies that $W_{d_{1},g_{1}}W_{d_{2},g_{2}}(X,g)=0$ if $g \neq g_{1}g_{2}$. Moreover the same equation implies that 
\begin{align*}
W_{d_{1},g_{1}}W_{d_{2},g_{2}}(X,g_1 g_2)&=1_{\Omega g_{1}^{-1}}(X)(d_{1}+I_{X})(\beta_{g_{1}}(d_{2})+I_{X})\\
&=1_{\Omega g_{1}^{-1}}d_{1}\beta_{g_{1}}(d_{2})+I_{X}.
\end{align*}
This proves $(3)$. 

Using $(1)$ and $(3)$, note that that for $d_{1},d_{2} \in \widetilde{D}$, 
\begin{align*}
W_{d_{1},e}W_{d_{2},e}&=W_{1_{\Omega}d_{1}d_{2},e} \\
                      &=W_{d_{1}d_{2},e}
\end{align*}
Thus the map $\widetilde{D} \ni d \to W_{d,e} \in \mathcal{D} \rtimes \clg$ is multiplicative. That it preserves the adjoint follows from $(2)$. Now let $d \in \widetilde{D}$. Note that  $W_{d,e}=0$ if and only if $1_{\Omega}d=0$ i.e. if and only if $d(a)=0$ for every $a \in P$. By Remark \ref{the kernel}, it follows that the kernel of the map $\widetilde{D} \ni d \to W_{d,e} \in \mathcal{D} \rtimes \clg$ is $I_{\Omega}$. 
This completes the proof. \hfill $\Box$.

In \cite{Li-Cuntz-Echterhoff} [Thm 5.2], it is proved that if $(P,G)$ is a quasi-lattice ordered pair and if the action $\alpha$ of $P$ on $A$ is obtained by restricting an action of $G$ on $A$ then the inclusion $A \ni x \to x \in A \rtimes_{red} P$ induces isomorphism at the $K$-theory level under certain assumptions on $G$.  We prove a very modest result that when $P=\mathbb{F}_{n}^{+}$, the free semigroup on $n$-generators, then the inclusion $A \ni x \to x \in A \rtimes_{red} P$ induces isomorphism at the $K$-theory level. This relies on the fact that $A \rtimes_{red} \mathbb{F}_{n}^{+}$ admits a presentation in terms of generators and relations which forms the content of the next section.
We restrict our attention only to the unital case. The isomorphism $K_{*}(A) \cong K_{*}(A \rtimes_{red} \mathbb{F}_{n}^{+})$ can be deduced from the unital case and by making use of Lemma \ref{unitisation}

\section{Quasi-lattice ordered case}
In this section, we show that the full groupoid crossed product $\mathcal{D} \rtimes \mathcal{G}$ admits a presentation in terms of generators and relations when $(P,G)$ is quasi-lattice ordered and the action $\alpha:P \to End(A)$ is unital. Throughout this section, we assume that $(P,G)$ is a quasi-lattice ordered pair, $A$ is a unital $C^{*}$-algebra and $\alpha:P \to End(A)$ is a unital action. 
\begin{dfn}
 Let $B$ be a unital $C^{*}$-algebra. Denote the isometries of $B$ by $\mathcal{V}(B)$. Let $\pi:A \to B$ be a  $*$-homomorphism and $V:P \to \mathcal{V}(B)$ be an anti-homomorphism i.e. $V_{e}=1$ and $V_{a}V_{b}=V_{ba}$ for $a,b \in P$. The pair $(\pi,V)$ is called Nica covariant if 
\item[(1)] for $x \in A$ and $a \in P$, $\pi(x)V_{a}=V_{a}\pi(\alpha_{a}(x))$,
\item[(2)] for $a,b \in P$, \begin{equation*}
 E_{a}E_{b}:=\begin{cases}
 E_{c} & \mbox{ if
} Pa \cap Pb = Pc,\cr
   &\cr
    0 &  \mbox{if } Pa \cap Pb =\emptyset.
         \end{cases}
\end{equation*} where $E_{a}:=V_{a}V_{a}^{*}$ for $a \in P$.
If $(\pi,V)$ is Nica covariant, we say $(\pi,V)$ is a Nica covariant representation of the triple $(A,P,\alpha)$. 
\end{dfn}

It is convenient to introduce the following universal algebra. 
\begin{dfn}
\label{defining relations}
 We let $A \rtimes P$ be the universal unital $C^{*}$-algebra generated by a unital copy of $A$ and isometries $\{v_{a}:a \in P\}$ such that 
\begin{enumerate}
\item[(C1)] for $x \in A$ and $a \in P$, $xv_{a}=v_{a}\alpha_{a}(x)$,
\item[(C2)] for $a,b \in P$, $v_{a}v_{b}=v_{ba}$, and
\item[(C3)] for $a, b \in P$, \begin{equation*}
 e_{a}e_{b}:=\begin{cases}
 e_{c} & \mbox{ if
} Pa \cap Pb = Pc,\cr
   &\cr
    0 &  \mbox{if } Pa \cap Pb =\emptyset.
         \end{cases}
\end{equation*} where $e_{a}:=v_{a}v_{a}^{*}$ for $a \in P$.
\end{enumerate}
\end{dfn}
Let $(\pi,V)$ be a Nica covariant representation of the triple $(A,P,\alpha)$ on a unital $C^{*}$-algebra $B$ and assume that $\pi$ is unital. Then there exists a unique $*$-homomorphism $\pi \rtimes V: A \rtimes P \to B$ such that $(\pi \rtimes V)(x)=\pi(x)$ and $(\pi \rtimes V)(v_{a})=V_{a}$ for $a \in P$. 
\begin{rmrk}
\label{regular}
Consider the pair $(\pi,V)$ constructed in Section 3 on the Hilbert module $A$-module $A \otimes \ell^{2}(P)$. Recall that for $x \in A$ and $a \in P$, the operators $\pi(x)$ and $V_{a}$ are given by:
\begin{align*}
\pi(x)(y \otimes \delta_{b})&:=\alpha_{b}(x)(y)\otimes \delta_{b} \\
V_{a}(y \otimes \delta_{b})&:=y \otimes \delta_{ba}
\end{align*}
We leave it to the reader to verify that $(\pi,V)$ is Nica covariant. As a consequence, it follows that $A \rtimes P$ is non-zero. In particular, we obtain a unital $*$-homomorphism $\rho:A \rtimes P \to A \rtimes_{red} P$ such that $\rho(x)=\pi(x)$ and $\rho(v_{a})=V_{a}$ for every $x \in A$ and $a \in P$. Let us call the map $\rho$ as the regular representation of $A \rtimes P$ and the pair $(\pi,V)$ the standard Nica-covariant pair of $(A,P,\alpha)$. 
\end{rmrk}

The following relation satisfied by the isometries $\{v_{a}:a \in P\}$ is due to Nica. We include the proof for completeness. Let $a,b \in P$ be given.  Then \begin{equation}
\label{isometries}
v_{a}^{*}v_{b}:=\begin{cases}
v_{ca^{-1}}v_{cb^{-1}}^{*} & \mbox{ if}
 Pa \cap Pb = Pc,\cr
  &\cr
  0 &  \mbox{if } Pa \cap Pb =\emptyset.
      \end{cases}
\end{equation} 
To see this observe that if $Pa \cap Pb=\emptyset$ then $e_{a}e_{b}=v_{a}v_{a}^{*}v_{b}v_{b}^{*}=0$ implies that $v_{a}^{*}v_{b}=0$. Now suppose $Pa \cap Pb =Pc$. Choose $s,t \in P$ such that $sa=tb=c$. Now calculate as follows to observe that 
\begin{align*}
v_{a}^{*}v_{b}&=v_{a}^{*}v_{a}v_{a}^{*}v_{b}v_{b}^{*}v_{b} \\
              &=v_{a}^{*}e_{a}e_{b}v_{b} \\
              &=v_{a}^{*}e_{c}v_{b} \\
              &=v_{a}^{*}v_{sa}v_{tb}^{*}v_{b} \\
              &=v_{a}^{*}v_{a}v_{s}v_{t}^{*}v_{b}^{*}v_{b} \\
              &=v_{ca^{-1}}v_{cb^{-1}}^{*}.
\end{align*}
Let us introduce some notation which will be used throughout this section. Consider the $C^{*}$-algebra $A \rtimes P$.  For $g \in G$, let 
\begin{equation}
w_{g}:= \begin{cases}
 0 & \mbox{if~} Pg \cap P=\emptyset, \cr
 & \cr 
 v_{a}v_{ag^{-1}}^{*} & \mbox{if~} Pg \cap P =Pa.
  \end{cases}
\end{equation}
Note that for $g \in G$, $w_{g}$ is a product of two partial isometries with commuting range projections. Hence  $w_{g}$ is a partial isometry
for each $g \in G$. Also note that $w_{g}^{*}=w_{g^{-1}}$. For $g \in G$, let $e_{g}:=w_{g}w_{g}^{*}$. Note that by definition if $w_{g} \neq 0$ then $e_{g}=e_{a}$ for some $a \in P$. Thus $\{e_{g}:g \in G\}$ forms a commuting family of projections.  The relations satisfied by $\{w_{g}:g \in G\}$ and $\{e_{g}: g \in G\}$ are summarised in the following lemma.

\begin{lmma}
\label{w-relations}
With the foregoing notations, 
\begin{enumerate}
\item[(1)] for $g, h \in G$, 
 \begin{equation*}
e_{g}e_{h}:=\begin{cases}
          0 & \mbox{if } Pg \cap Ph \cap P =\emptyset, \cr
            & \cr 
            e_{c} & \mbox{if }Pg \cap Ph \cap P=Pc ,
            \end{cases}
            \end{equation*}
\item[(2)] for $g, h \in G$, $w_{g}w_{h}=e_{g}w_{hg}$.             
\end{enumerate}
\end{lmma}
\textit{Proof.} Let $g,h \in G$ be given. Suppose  $Pg \cap Ph \cap P\neq \emptyset$. Then $Pg \cap P \neq \emptyset$ and $Ph \cap P \neq \emptyset$. Choose $a,b,c \in P$ be such that $Pg \cap P=Pa$ and $Ph \cap P=Pb$ and $Pg \cap Ph \cap P=Pa \cap Pb = Pc$. By definition, $e_{g}=e_{a}$ and $e_{h}=e_{b}$. Thus $e_{g}e_{h}=e_{c}$.  Now suppose $Pg \cap Ph \cap P=\emptyset$. We claim that $e_{g}e_{h}=0$. Suppose $e_{g}e_{h} \neq 0$. Then $e_{g} \neq 0$ and $e_{h} \neq 0$. Thus $Pg \cap P \neq \emptyset$ and $Ph \cap P \neq \emptyset$. Choose $a,b \in P$ such that $Pg \cap P=Pa$ and $Ph \cap P=Pb$. Note that by definition $e_{g}=e_{a}$ and $e_{h}=e_{b}$. Now $e_{g}e_{h}=e_{a}e_{b} \neq 0$ implies that $Pa \cap Pb=Pg \cap P \cap Ph$ is non-empty which is a contradiction. Thus if $Pg \cap Ph \cap P=\emptyset$ then $e_{g}e_{h}=0$. This proves $(1)$. 

Let $g, h \in G$ be given. \textit{Claim:} $w_{g}w_{h}=0 \Leftrightarrow Phg \cap Pg \cap P =\emptyset \Leftrightarrow e_{g}w_{gh}=0$. Note the following chain of equivalences.
\begin{align*}
w_{g}w_{h}=0 & \Leftrightarrow w_{g}w_{h}w_{h}^{*}w_{g}^{*}=w_{g}e_{h}e_{h}w_{g}^{*}=0 \\
             & \Leftrightarrow e_{h}w_{g^{-1}}=0 \\
             & \Leftrightarrow e_{h}w_{g^{-1}}w_{g^{-1}}^{*}e_{h}=0 \\
             & \Leftrightarrow e_{h}e_{g^{-1}}e_{h}=e_{h}e_{g^{-1}}=0 \\
             & \Leftrightarrow Ph \cap Pg^{-1} \cap P = \emptyset ~~(\textrm{by (1)}) \\
             & \Leftrightarrow Phg \cap P \cap Pg = \emptyset \\
             & \Leftrightarrow e_{g}e_{hg}=e_{g}e_{hg}e_{g}=e_{g}w_{hg}w_{hg}^{*}e_{g}=0 ~~(\textrm{ by (1)}) \\ 
             & \Leftrightarrow e_{g}w_{hg}=0.
\end{align*}
This proves the claim. Thus to prove (2), we can and will  assume that $w_{g}w_{h} \neq 0$ and $e_{g}w_{hg} \neq 0$. This in particular implies that $Pg \cap P$,  $Ph \cap P$, $Phg  \cap P$ and  $Phg \cap Pg \cap P$ are non-empty. Choose $a,b,c, d \in P$ such that $Pg \cap P=Pa$,  $Ph \cap P=Pb$, $Phg \cap Pg \cap P=Pc$ and $Phg \cap P=Pd$. 

Let $r,s,t \in P$ be such that $rg=a$, $sh=b$ and $thg=d$. Now note that $Pr \cap Pb = Pag^{-1} \cap P \cap Ph=Pg^{-1} \cap P \cap Ph=Pcg^{-1}$. Hence $cg^{-1} \in P$. Also note that $Pa \cap Pd=Pg \cap P \cap Phg \cap P=Phg \cap Pg \cap P=Pc$.

Now compute as follows to observe that 
\begin{align*}
w_{g}w_{h}&= v_{a}v_{r}^{*}v_{b}v_{s}^{*} \\
          &=v_{a}v_{cg^{-1}r^{-1}}v_{cg^{-1}b^{-1}}^{*}v_{s}^{*} ~~(\textrm{by Eq. \ref{isometries}}) \\
          &=v_{cg^{-1}r^{-1}a}v_{cg^{-1}b^{-1}s}^{*} \\
          &=v_{ca^{-1}a}v_{cg^{-1}h^{-1}}^{*} \\
          &=v_{c}v_{c(hg)^{-1}}^{*}
\end{align*}
and 
\begin{align*}
e_{g}w_{gh}&=e_{a}v_{d}v_{t}^{*} \\
           &=v_{a}v_{a}^{*}v_{d}v_{t}^{*} \\
           &=v_{a}v_{ca^{-1}}v_{cd^{-1}}^{*}v_{t}^{*} ~~(\textrm{by Eq. \ref{isometries}})\\
           &=v_{c}v_{cd^{-1}t}^{*} \\
           &=v_{c}v_{c(hg)^{-1}}^{*}.
\end{align*}
Hence $w_{g}w_{h}=e_{g}w_{hg}$. This proves $(2)$. This completes the proof. \hfill $\Box$.

Let $(\mathcal{D},\clg,\beta)$ be the Wiener-Hopf groupoid dynamical system associated to  $(A,P,G,\alpha)$. Since we are dealing with the unital case, by Remark \ref{unital case}, the dynamical system $(\mathcal{D},\clg,\beta)$ is given as follows. For $g \in G$ and $x \in A$, let $j_{g}(x) \in \ell^{\infty}(G,A)$ be defined by \begin{equation*}
 j_{g}(x)(h):=\begin{cases}
 \alpha_{hg^{-1}}(x) & \mbox{ if
} h \in Pg,\cr
   &\cr
    0 &  \mbox{ if } h \notin Pg.
         \end{cases}
\end{equation*}
Let $\widetilde{D}$ be the $C^{*}$-algebra generated by $\{j_{g}(x): g \in G, x \in A\}$. The  $C^{*}$-algebra $\widetilde{D}$ is invariant under the translation action $\beta$ and contains $C_{0}(\widetilde{\Omega})$ as a sub $C^{*}$-algebra where we identify $C_{0}(\widetilde{\Omega})$ with the $C^{*}$-subalgebra generated by $\{j_{g}(1):g \in G\}$. Thus $\widetilde{D}$ is a $(\widetilde{\Omega},G)$ algebra. 

 Let $\displaystyle \mathcal{D}:=\coprod_{X \in \Omega} \widetilde{D}/I_{X}$ be the corresponding upper semicontinuous bundle over $\Omega$ and let   $(\beta_{(X,g)})_{(X,g) \in \clg}$ be the action of the Wiener-Hopf groupoid $\clg:= \widetilde{\Omega} \rtimes G|_{\Omega}$ which  is given by $\beta_{(X,g)}:\widetilde{D}/I_{Xg} \ni d+I_{Xg} \to \beta_{g}(d)+I_{X} \in \widetilde{D}/I_{X}$.  

We claim that the linear span of $\{j_{g}(x):g \in G,x \in A\}$ is a dense $*$-subalgebra of $\widetilde{D}$. To see this, let $g_{1},g_{2} \in G$ and $x_{1},x_{2} \in A$ be such that $j_{g_{1}}(x_{1})j_{g_{2}}(x_{2}) \neq 0$. Then $Pg_1 \cap Pg_2 \neq \emptyset$. Since $(P,G)$ is quasi-lattice ordered, it follows that there exists $g \in G$ such that $Pg_{1} \cap Pg_{2}=Pg$. Choose $a,b \in P$ such that $ag_{1}=g$ and $bg_{2}=g$. Then $j_{g_{1}}(x_{1})j_{g_{2}}(x_{2})=j_{g}(\alpha_{a}(x_{1})\alpha_{b}(x_{2}))$.

\begin{rmrk}
Note that by Lemma \ref{Relations}, it follows that for $d \in \widetilde{D}$ and $g \in G$,
 \begin{align*}
 W_{1_{\Omega},e}W_{d,g}&=W_{1_{\Omega}d,g}~~ (\textrm{ by $(3)$, Lemma \ref{Relations})} \\
 &=W_{d,g}~~ (~~\textrm{by $(1)$, Lemma \ref{Relations})}
 \end{align*}
  Thus $W_{1_{\Omega},e}$ is the multiplicative identity of $\Gamma_{c}(\clg,r^{*}\mathcal{D})$.
  \end{rmrk}
 \begin{ppsn}
 \label{homomorphism lambda}
 With the foregoing notations, there exists a  $*$-homomorphism $\lambda:A \rtimes P \to \mathcal{D} \rtimes \clg$ such that $\lambda(v_{a})=W_{1_{\Omega},a^{-1}}$ and $\lambda(x)=W_{j_{e}(x),e}$ for $a \in P$ and $x \in A$. Moreover $\lambda$ is onto.
  \end{ppsn}
\textit{Proof.} Let $\pi:A \to \Gamma_{c}(\clg,r^{*}\mathcal{D})$ be defined by $\pi(x)=W_{j_{e}(x),e}$. By Lemma \ref{Relations} (Part (3)), it follows that $\pi$ is a $*$-representation. Also note that $\pi$ is unital. For $a \in P$, let $\widetilde{v}_{a}:=W_{1_{\Omega},a^{-1}}$ and let $\widetilde{e}_{a}:=\widetilde{v}_{a}\widetilde{v}_{a}^{*}$.  We verify that $(\pi,\widetilde{v})$ is a Nica covariant. Let $a \in P$ be given. By Lemma \ref{Relations}, it follows that 
\begin{align*}
\widetilde{v}_{a}^{*}\widetilde{v}_{a}&=W_{\beta_{a}(1_{\Omega}),a}W_{1_{\Omega},a^{-1}} \\
                                      &=W_{1_{\Omega a^{-1}},a}W_{1_{\Omega},a^{-1}} \\
                                     &=W_{1_{\Omega a^{-1}}\beta_{a}(1_{\Omega}),e} \\
                                     &=W_{1_{\Omega a^{-1}},e} \\
                                     &=W_{1_{\Omega},e} ~~(\textrm{by (1) of Lemma \ref{Relations}})   
  \end{align*}
Thus $\{\widetilde{v}_{a}: a \in P\}$ is a collection of isometries. The fact that $(\pi,\widetilde{v})$ is Nica covariant follows from  repeated applications of Lemma \ref{Relations}.  

Let $x \in A$ and $a \in P$ be given. First by definition, $j_{a}(\alpha_{a}(x))=j_{e}(x)1_{\Omega}1_{\Omega a}$. Now note  that 
\begin{align*}
\pi(x)\widetilde{v}_{a}&=W_{j_{e}(x),e}W_{1_{\Omega},a^{-1}} \\
                       &=W_{j_{e}(x)1_{\Omega},a^{-1}} \\
                       &=W_{j_{e}(x)1_{\Omega}1_{\Omega a},a^{-1}} \\
                       &=W_{j_{a}(\alpha_{a}(x)),a^{-1}}\\
                       &=W_{1_{\Omega}1_{\Omega a}j_{a}(\alpha_{a}(x)),a^{-1}} \\
                       &=W_{1_{\Omega}1_{\Omega a}\beta_{a}^{-1}(j_{e}(\alpha_{a}(x)),a^{-1}}\\
                       &=W_{1_{\Omega},a^{-1}}W_{j_{e}(\alpha_{a}(x)),e} \\
                       &=\widetilde{v}_{a}\pi(\alpha_{a}(x)).
\end{align*}  
 We leave it to the reader to verify that  $\widetilde{v}_{a}\widetilde{v}_{b}=\widetilde{v}_{ba}$ for $a,b \in P$. Another application of Lemma \ref{Relations} imply that $\widetilde{e}_{a}:=\widetilde{v}_{a}\widetilde{v}_{a}^{*}=W_{1_{\Omega a},e}=W_{j_{a}(1),e}$ for every $a$. Now note that for $a,b \in P$, $j_{a}(1)j_{b}(1)=0$ if $Pa \cap Pb = \emptyset$ and $j_{a}(1)j_{b}(1)=j_{c}(1)$ if $Pa \cap Pb=Pc$. Since the map $\widetilde{D} \ni d \to W_{d,e} \in \mathcal{D} \rtimes \clg$ is a $*$-homomorphism, it follows that 
\begin{equation*}
 \widetilde{e}_{a}\widetilde{e}_{b}:=\begin{cases}
 \widetilde{e}_{c} & \mbox{ if
} Pa \cap Pb = Pc,\cr
   &\cr
    0 &  \mbox{if } Pa \cap Pb =\emptyset.
         \end{cases}
\end{equation*} 
 Thus by the universal property of $A \rtimes P$, there exists a unique $*$-homomorphism $\lambda:A \rtimes P \to \mathcal{D} \rtimes \clg$ such that $\lambda(x)=W_{j_{e}(x),e}$ for $x \in A$ and $\lambda(v_{a})=W_{1_{\Omega},a^{-1}}$ for $a \in P$.
 
Let $\mathcal{C}$ be the $C^{*}$-subalgebra of $\mathcal{D} \rtimes \clg$ generated by $\{W_{j_{e}(x),e}:x \in A\}$ and $\{\widetilde{v}_{a}:a \in P\}$.  To show that $\lambda$ is surjective, it is enough to show that for each $d \in \widetilde{D}$ and $g \in G$, $W_{d,g} \in \mathcal{C}$. 

\textit{Claim:} $W_{d,e} \in \mathcal{C}$ for every $d \in \widetilde{D}$. Since the map $\widetilde{D} \ni d \to W_{d,e} \in \mathcal{D} \rtimes \clg$ is a $*$-homomorphism, it is enough to show that $W_{j_{g}(x),e}=W_{1_{\Omega}j_{g}(x),e} \in \mathcal{C}$ for every $g \in G$ and $x \in A$. Let $x \in A$ and $g \in G$ be such that $1_{\Omega}j_{g}(x) \neq 0$. Then, by definition, $Pg \cap P \neq \emptyset$. Since $(P,G)$ is a quasi-lattice ordered pair, it follows that there exists $a \in P$ such that $Pg \cap P = Pa$.  Let $b \in P$ be such that $bg=a$. Let $h \in G$. If $h \notin Pa=Pg \cap P$ then $1_{\Omega}(h)j_{g}(x)(h)=0$. Suppose $h \in Pa$. Then \[
1_{\Omega}(h)j_{g}(x)(h)=\alpha_{hg^{-1}}(x)=\alpha_{ha^{-1}b}(x)=\alpha_{ha^{-1}}(\alpha_{b}(x))=j_{a}(\alpha_{b}(x)).\] Hence $1_{\Omega}j_{g}(x)=j_{a}(\alpha_{b}(x))$. Thus to prove the claim it is enough to show that $W_{j_{a}(x),e} \in \mathcal{C}$ for every $a \in P$ and $x \in A$. But observe using Lemma \ref{Relations} that $W_{j_{a}(x),e}=\widetilde{v}_{a}W_{j_{e}(x),e}\widetilde{v}_{a}^{*}$ for $x \in A$ and $a \in P$. This proves the claim.
 
Let $g \in G$ be given. Observe that by $(1)$ of Lemma \ref{Relations}, $W_{1_{\Omega},g} \neq 0$ if and only if $1_{\Omega}1_{\Omega g^{-1}} \neq 0$ i.e. if and only if $Pg^{-1} \cap P \neq \emptyset$. Let $a \in P$ be such that $Pg^{-1}\cap P=Pa$ and choose $b \in P$ such that $bg^{-1}=a$. Observe that $1_{\Omega}(h)1_{\Omega g^{-1}}(h)=1$ if and only if $h \in Pg^{-1} \cap P=Pa$. Thus $1_{\Omega}1_{\Omega g^{-1}}=1_{\Omega a}$. 

Again by Lemma \ref{Relations}, note that 
\begin{align*}
\widetilde{v}_{a}\widetilde{v}_{b}^{*}&=W_{1_{\Omega},a^{-1}}W_{1_{\Omega b^{-1}},b} \\
                                      &=W_{1_{\Omega}1_{\Omega a}1_{\Omega b^{-1}a},a^{-1}b} \\
                                      &= W_{1_{\Omega a},a^{-1}b} \\
                                      &=W_{1_{\Omega}1_{\Omega g^{-1}},g} \\
                                      &=W_{1_{\Omega},g}~~(\textrm{by (1) of Lemma \ref{Relations}}).
\end{align*}
This implies that $\{W_{1_{\Omega},g}:g \in G\}$ is contained in $\mathcal{C}$. Now let $d \in \widetilde{D}$ and $g \in G$ be given. Note that $W_{d,g}=W_{d,e}W_{1_{\Omega},g}$. As a consequence, it follows that $\{W_{d,g}:d \in \widetilde{D}, g \in G\}$ is contained in $\mathcal{C}$. Since the linear span of $\{W_{d,g}:d \in \widetilde{D}, g \in G\}=\Gamma_{c}(\clg,r^{*}\mathcal{D})$ is dense in $\mathcal{D} \rtimes \clg$, it follows that $\mathcal{C}=\mathcal{D} \rtimes \clg$. This proves that $\lambda$ is surjective. This completes the proof. \hfill $\Box$.

The main aim of this section is to show that the map $\lambda$ of the preceeding proposition is an isomorphism. We do this by constructing the inverse of $\lambda$.

For $g \in P$ and $x \in A$, let $i_{g}(x) \in \ell^{\infty}(P,A)$ be the restriction of $j_{g}(x)$ onto $P$. Let $g \in G$ and $x \in A$ be such that $i_{g}(x) \neq 0$. Then $Pg \cap P \neq \emptyset$. Let $a \in P$ be such that $Pg \cap P=Pa$ and let $b \in P$ be such that $bg=a$. Note that $i_{g}(x)=i_{a}(\alpha_{b}(x))$. Thus the $C^{*}$-algebra generated by $\{i_{g}(x): g \in G, x \in A\}$ coincides with the $C^{*}$-subalgebra generated by $\{i_{a}(x):x \in A, a \in P\}$.  Let $D$ denote the $C^{*}$-subalgebra generated by $\{i_{a}(x):a \in P, x \in A\}$. Note by Remark \ref{the kernel} and $(4)$ of Lemma \ref{Relations}, the map, call it $\sigma$, $\Gamma(\Omega,\mathcal{D}) \ni W_{d,e} \to (d(a))_{a \in P} \in \ell^{\infty}(P,A)$ is injective with range $D$.

 Let $\mathcal{E}$ be the $C^{*}$-subalgebra of $A \rtimes P$ generated by the projections $\{e_{a}:a \in P\}$ and let $\mathcal{F}$ be the $C^{*}$-subalgebra of $A \rtimes P$ generated by $\{v_{a}xv_{a}^{*}:x \in A, a \in P\}$. We first prove that $\lambda$ is $1$-$1$ on $\mathcal{F}$. Let $\widetilde{\lambda}:\mathcal{F} \to D$ be defined by $\widetilde{\lambda}=\sigma \circ \lambda$. Note that for $a \in P$ and $x \in A$, $\lambda(v_{a}xv_{a}^{*})=W_{j_{a}(x),e}$ and $\widetilde{\lambda}(v_{a}xv_{a}^{*})=i_{a}(x)$.  We show $\lambda$ is injective on $\mathcal{F}$ by showing that $\widetilde{\lambda}$ is injective.  
 
 \begin{rmrk}
 Observe the following. 
 \begin{enumerate}
 \item[(1)] For $g \in G$, $e_{g} \in \mathcal{E}$. For if $e_{g} \neq 0$ then $e_{g}=e_{a}$ where $a \in P$ is such that $Pg \cap P=Pa$. 
 \item[(2)] For $g \in G$ and $x \in A$, $w_{g}xw_{g}^{*} \in \mathcal{F}$. For, suppose $x \in A$ and $g \in G$ be such that $w_{g} \neq 0$. Choose $a,b \in P$ such that $Pg \cap P=Pa$ and $bg=a$. Then by definition $w_{g}xw_{g}^{*}=v_{a}v_{b}^{*}xv_{b}v_{a}^{*}=v_{a}\alpha_{b}(x)v_{a}^{*} \in \mathcal{F}$. 
 \item[(3)] We leave it to the reader to verify that $\lambda(w_{g}xw_{g}^{*})=W_{j_{g}(x),e}$ and $\widetilde{\lambda}(w_{g}xw_{g}^{*})=i_{g}(x)$ for $g \in G$ and $x \in A$. 
 \end{enumerate}
 
 \end{rmrk}

 The spectrum of $\mathcal{E}$ can be identified with the Wiener-Hopf compactification $\Omega$ of $(P,G)$. This is essentially due to Nica (give references). We recall here the proof for completeness in the next few paragraphs. Notation: As $A$ is unital, we identify $\ell^{\infty}(P)$ as a $*$-subalgebra of $\ell^{\infty}(P,A)$. For a subset $Y \subset P$, let $1_{Y} \in \ell^{\infty}(P)$ denote the characteristic function of $Y$. Note that $\widetilde{\lambda}(e_{a})=1_{Pa}$ for $a \in P$. We also identify $C(\Omega)$ as a $*$-subalgebra of $\ell^{\infty}(P)$ via the map $C(\Omega) \ni \phi \to (\phi(P^{-1}a))_{a \in P} \in \ell^{\infty}(P)$. Under this map $1_{\Omega a}$ is mapped to $1_{Pa}$ for $a \in P$.
 
 \begin{enumerate}
 \item[(1)] For $a, b \in P$, let $a \leq b$ if there exists $c \in P$ such that $b=ca$. Let $A$ be a subset of $P$. The subset $A$ is called \textbf{directed} if given $a,b \in A$ there exists $c \in A$ such that $c \geq a,b$. $A$ is called \textbf{hereditary} if $b \in A$ and $a \leq b$ then $a \in A$. 
 \item [(2)] Let $\Omega_{\mathcal{N}}$ denote the set of all non-empty, directed and hereditary subsets of $P$. We topologise $\Omega_{\mathcal{N}}$ by considering $\Omega_{\mathcal{N}}$ as a subset of $\{0,1\}^{P}$ and endow $\Omega_{\mathcal{N}}$ with the subspace topology inherited from the product topology on $\{0,1\}^{P}$.  For $a \in P$, let $[e,a]:=\{x \in P: x \leq a\}$. Then $\{[e,a]: a \in P\}$ is dense in $\Omega_{\mathcal{N}}$. For $a \in P$, let $U_{a}:=\{A \in \Omega_{\mathcal{N}}: a \in A\}$. 
  \item[(3)] Let $\Omega$ be the Wiener-Hopf compactification of $(P,G)$. Let $A \in \Omega$ be given. Note that  $A$ contains the identity element $e$ and $P^{-1}A \subset A$. Thus $A \cap P$ is non-empty and hereditary. We claim that $A \cap P$ is directed. Let $a,b \in A \cap P$ be given. Choose a sequence $(a_{n})$ in $P$ such that $P^{-1}a_{n} \to A$. Then there exists $N \in \mathbb{N}$ such that $a,b \in P^{-1}a_{n}$ for $n \geq N$. In other words, $a_{n} \in Pa \cap Pb$ eventually. Since $(P,G)$ is quasi-lattice ordered, it follows that there exists $c \in P$ such that $Pa \cap Pb = Pc$. Thus $a_{n} \in Pc$ eventually or $c \in P^{-1}a_{n}$ eventually. Thus $c \geq a,b $ and $c \in A \cap P$. This proves that $A \cap P$ is non-empty. Thus the map $\Omega \ni A \to A \cap P \in \Omega_{\mathcal{N}}$ is well-defined.
  \item[(4)] The map $\Omega \ni A \to A \cap P \in \Omega_{\mathcal{N}}$ is a homeomorphism. Clearly the map $\Omega \ni A \to A \cap P \in \Omega_{\mathcal{N}}$ is continuous. Let $A,B \in \Omega$ be such that $A \cap P=B \cap P$. Choose a sequence $a_{n} \in P$ such that $P^{-1}a_{n} \to A$. Consider an element $g \in A$. Then $a_{n} \in Pg$ eventually. Thus $Pg \cap P \neq \emptyset$. Since $(P,G)$ is quasi-lattice ordered, it follows that there exists $a \in P$ such that $Pg \cap P=Pa$. Then $a_{n} \in Pa$ eventually. Thus $a \in A \cap P=B \cap P$. This implies that $a \in B$. Note that $g \leq a$ and $P^{-1}B \subset B$. Thus $g \in B$. This shows that $A \subset B$. Similarly $B \subset A$. This proves that the map $\Omega \ni A \to A \cap P$ is $1-1$. Note that $P^{-1}a \cap P=[e,a]$ for $a \in P$. Since $\{[e,a]: a \in P\}$ is dense in $\Omega_{\mathcal{N}}$ and $\Omega$ is compact, it follows that the map $\Omega \ni A \to A \cap P \in \Omega_{\mathcal{N}}$ is a homeomorphism and we identify $\Omega$ and $\Omega_{\mathcal{N}}$ via this homeomorphism. Under the homeomorphism $\Omega \ni A \to A \cap P \in \Omega_{\mathcal{N}}$, by Lemma \ref{facts on Omega}, the subset $\Omega a$ is mapped onto $U_{a}$ for every $a \in P$. 
  
  \item[(5)] Let $\widehat{\mathcal{E}}$ be the character space of $\mathcal{E}$. For $\chi \in \widehat{\mathcal{E}}$, let $A_{\chi}:=\{a \in P: \chi(e_{a})=1\}.$ Then $A_{\chi}$ contains the identity element and $A_{\chi}$ is directed and hereditary. The map $\widehat{\mathcal{E}} \ni \chi \to A_{\chi} \in \Omega_{\mathcal{N}} \cong \Omega$ is clearly continuous and injective. Thus one obtains a $*$-homomorphism $\widetilde{\mu}: C(\Omega) \to C(\widehat{\mathcal{E}}) \cong \mathcal{E}$ defined by $\widetilde{\mu}(f)(\chi)=f(B_{\chi})$ for $f \in C(\Omega)$ and $\chi \in \widehat{\mathcal{E}}$ where $B_{\chi}$ is the unique element in $\Omega$ such that $B_{\chi}\cap P = A_{\chi}$. We leave it to the reader to verify that $\widetilde{\mu}(1_{\Omega a})=e_{a}$ for $a \in P$. Since $\widetilde{\lambda}(e_{a})=1_{Pa}$, it follows that $\widetilde{\lambda}: \mathcal{E} \to C(\Omega)$ is an isomorphism with $\widetilde{\mu}$ being its inverse. 
\end{enumerate} 
Let $\mathcal{I}:=\{Y \subset P: 1_{Y} \in C(\Omega) \subset \ell^{\infty}(P)\}$. The collection $\mathcal{I}$ is called the set of constructible left ideals of $P$ by Li in \cite{Li13}. 
Note that $\mathcal{I}$ is closed under finite unions, finite intersections and complements. Moreover for every $a \in P$, $Pa \in \mathcal{I}$. For $Y \in \mathcal{I}$, let $e_{Y} \in \mathcal{E}$ be such that $\widetilde{\lambda}(e_{Y})=1_{Y}$. For $Y \in \mathcal{I}$, note that $e_{Y}$ is a projection. 

\begin{lmma}
\label{structure of mathcal{F}}
Let $\widetilde{\mathcal{F}}$ be the linear span of  $\{v_{a}xv_{a}^{*}: x \in A, a \in P\}$. 
\begin{enumerate}
\item[(1)]  Then $\widetilde{\mathcal{F}}$ is a dense $*$-subalgebra of $\mathcal{F}$.
\item[(2)] The algebra $\mathcal{E}$ is contained in the center of $\mathcal{F}$.
\item[(3)] Let $z \in \widetilde{\mathcal{F}}$ be given. Then there exists $n \in \mathbb{N}$, $Y_{1},Y_{2},\cdots Y_{n} \in \mathcal{I}$, $a_{1},a_{2},\cdots a_{n} \in P$ and $x_{1},x_{2},\cdots x_{n} \in A$ such that $a_{i} \in Y_{i}$ for every $i$, $Y_{i} \cap Y_{j} =\emptyset$ if $i \neq j$ and $z=\sum_{i=1}^{n}e_{Y_{i}}v_{a_{i}}x_{i}v_{a_{i}}^{*}$.
\item[(4)] The map $\widetilde{\lambda}:\mathcal{F} \to D$ is isometric and onto.  
\end{enumerate}
\end{lmma}
\textit{Proof.} Let $a, b \in P$ and $x,y \in A$ be given. Suppose  $Pa \cap Pb=\emptyset$. The relation $e_{a}e_{b}=0$ implies that $v_{a}^{*}v_{b}=0$. Hence $(v_{a}xv_{a}^{*})(v_{b}yv_{b}^{*})=0$. Now suppose $Pa \cap Pb \neq \emptyset$. Since $(P,G)$ is quasi-lattice ordered, it follows that there exists $c \in P$ such that $Pa \cap Pb = Pc$. Choose $s,t \in P$ such that $sa=tb=c$. Now calculate as follows to find that 
\begin{align*}
(v_{a}xv_{a}^{*})(v_{b}yv_{b}^{*})&=v_{a}xv_{s}v_{t}^{*}yv_{b}^{*} ~~(\textrm{by Eq. \ref{isometries}}) \\
                                  &=v_{a}v_{s}\alpha_{s}(x)\alpha_{t}(y)v_{t}^{*}v_{b}^{*} \\
                                  &=v_{sa}\alpha_{s}(x)\alpha_{t}(y)v_{tb}^{*}\\
                                  &=v_{c}\alpha_{s}(x)\alpha_{t}(y)v_{c}^{*}.
\end{align*}
This proves that $\widetilde{\mathcal{F}}$ is closed under multiplication. Clearly $\widetilde{\mathcal{F}}$ is $*$-closed. Thus $\widetilde{\mathcal{F}}$ is dense in $\mathcal{F}$. This proves $(1)$. 

Let $a,b \in P$ and $x \in A$ be given. The calculation done to prove $(1)$ implies that 
\begin{equation}
\label{relation}
 e_{a}v_{b}xv_{b}^{*}:=\begin{cases}
 v_{c}\alpha_{cb^{-1}}(x)v_{c}^{*} & \mbox{ if
} Pa \cap Pb = Pc,\cr
   &\cr
    0 &  \mbox{if } Pa \cap Pb =\emptyset
         \end{cases}
\end{equation} 
and 
\begin{equation*}
 v_{b}xv_{b}^{*}e_{a}:=\begin{cases}
 v_{c}\alpha_{cb^{-1}}(x)v_{c}^{*} & \mbox{ if
} Pb \cap Pa = Pc,\cr
   &\cr
    0 &  \mbox{if } Pb \cap Pa =\emptyset.
         \end{cases}
\end{equation*} 
Hence $e_{a}v_{b}xv_{b}^{*}=v_{b}xv_{b}^{*}e_{a}$ for every $a,b \in P$ and $x \in A$. This proves $(2)$. 

Let $z:=\sum_{i=1}^{n}v_{b_{i}}y_{i}v_{b_{i}}^{*} \in \widetilde{\mathcal{F}}$ be given with $b_{i} \in P$ and $y_{i} \in A$. For $B \subset \{1,2,\cdots n\}$, let $Z_{B}:=\bigcap_{i \in B}Pb_{i}$. Since $(P,G)$ is quasi-lattice ordered, it follows that if $Z_{B}$ is non-empty, then there exists $a_{B} \in P$ such that $Z_{B}=Pa_{B}$. For $B \subset \{1,2,\cdots,n\}$, let $Y_{B}:=Z_{B} \backslash \displaystyle \bigcup_{B \subsetneq C}Z_{C}$. We use here the convention that empty intersections is the whole set, in this case, $P$ and empty unions is empty. 
We leave it to the reader to verify that  $\{Y_{B}: B \subset \{1,2,\cdots n\}\}$ forms a disjoint collection of subsets of $P$ and for $B \subset \{1,2,\cdots,n\}$, $Z_{B}=\displaystyle \bigcup_{C \supset B}Y_{C}$. Also observe that if $Y_{B} \neq \emptyset$ then $a_{B} \in Y_{B}$. To see this, let $B \subset \{1,2,\cdots,n\}$ be such that $Y_{B} \neq \emptyset$. Note that $Y_{B} \subset Z_{B}=Pa_{B}$. Thus if $a_{B} \notin Y_{B}$ then there exists $C \supsetneq B$ such that $a_{B} \in Z_{C}=Pa_{C}$. But $Z_{C}$ is closed under left multiplication by $P$. Thus $Z_{B} \subset Z_{C}$ and  hence $Y_{B} =\emptyset$. Thus if $Y_{B} \neq \emptyset$ then $a_{B} \in Y_{B}$. Let $\mathcal{V}:=\{B \subset \{1,2,\cdots,n\}: Y_{B} \neq \emptyset\}$.  
 
Now calculate as follows to see that
\begin{align*}
z&=\sum_{i=1}^{n}e_{Pb_{i}}v_{b_{i}}y_{i}v_{b_{i}}^{*} \\
 & = \sum_{i=1}^{n} \displaystyle \big( \sum_{ B \in \mathcal{V}, i \in B}e_{Y_{B}} \big)v_{b_{i}}y_{i}v_{b_{i}}^{*} \\
 &= \sum_{B \in \mathcal{V}}e_{Y_{B}}e_{Z_{B}}\big (\sum_{i \in B}e_{Z_{B}}v_{b_{i}}y_{i}v_{b_{i}}^{*} \big) \\
 &= \sum_{B \in \mathcal{V}}e_{Y_{B}} \big (\sum_{i \in B}e_{a_{B}}v_{b_{i}}y_{i}v_{b_{i}}^{*} \big) \\
 &= \sum_{B \in \mathcal{V}}e_{Y_{B}} \big (\sum_{i \in B}v_{a_{B}}\alpha_{a_{B}a_{i}^{-1}}(y_{i})v_{a_{B}}^{*} \big)~~~(\textrm{ By Eq. \ref{relation} and $Pa_{B} \subset Pa_{i}$ if $i \in B$}) \\
 &= \sum_{B \in \mathcal{V}}e_{Y_{B}}v_{a_{B}}x_{B}v_{a_{B}}^{*} 
\end{align*}
where for $B \in \mathcal{V}$, $\displaystyle x_{B}:=\sum_{i \in B}\alpha_{a_{B}a_{i}^{-1}}(y_{i})$. This proves $(3)$.

To prove $(4)$, it is enough to show that  $\widetilde{\lambda}$ is isometric on $\widetilde{\mathcal{F}}$. Let $z \in \widetilde{\mathcal{F}}$ be given. By $(3)$ we can write $z=\sum_{i=1}^{n}e_{Y_{i}}v_{a_{i}}x_{i}v_{a_{i}}^{*}$ with $Y_{i}$ disjoint, $a_{i} \in Y_{i}$ and $x_{i} \in A$. Since $\{e_{Y_{i}}: i=1,2\cdots,n\}$ are orthogonal central projections in $\mathcal{F}$, it follows that \[||z|| \leq \displaystyle \max_{i \in \{1,2,\cdots n \}}||v_{a_{i}}x_{i}v_{a_{i}}^{*}|| \leq \displaystyle \max_{i \in \{1,2,\cdots,n\}}||x_{i}||.\] 

Note that $\widetilde{\lambda}(z)=\sum_{i=1}^{n}1_{Y_{i}}i_{a_{i}}(x_{i})$. Since  $a_i \in Y_i$,  $\widetilde{\lambda}(z)(a_{i})=x_{i}$ for $i \in \{1,2,\cdots,n\}$. Thus $||x_{i}|| \leq ||\widetilde{\lambda}(z)||$ for every $i \in \{1,2,\cdots,n\}$. As a consequence, it follows that $||z|| \leq || \widetilde{\lambda}(z)||$. But $\widetilde{\lambda}:\mathcal{F} \to D$ is a $*$-homomorphism and hence $||\widetilde{\lambda}(z)|| \leq ||z||$. This proves that $||\widetilde{\lambda}(z)||=||z||$ for every $z \in \widetilde{\mathcal{F}}$. Since $\widetilde{\mathcal{F}}$ is dense in $\mathcal{F}$, it follows that $\widetilde{\lambda}$ is isometric. Surjectivity of $\widetilde{\lambda}$ follows from the fact that $\widetilde{\lambda}(v_{a}xv_{a}^{*})=i_{a}(x)$ for $x \in A$ and $a \in P$ and the observation that $D$ is generated by $\{i_{a}(x):a \in P, x \in A\}$.  This completes the proof. \hfill $\Box$ 

Let $\widetilde{\mu}:D \to \mathcal{F}$ be the inverse of $\widetilde{\lambda}:\mathcal{F} \to D$. Let us denote the map $\widetilde{D} \ni d \to (d(a))_{a \in P} \in D$ by $res$ and let $\nu:\widetilde{D} \to \mathcal{F}$ be defined by $\nu=\widetilde{\mu} \circ res$. Note that $\nu(1_{\Omega g})=e_{g}$ for every $g \in G$. In particular, $\nu(1_{\Omega})=1$.

Now we define a map $\mu: \Gamma_{c}(\clg,r^{*}\mathcal{D}) \to A \rtimes P$ as follows: Let $f \in \Gamma_{c}(\clg,r^{*}\mathcal{D})$. Write $f=\displaystyle \sum_{d \in \widetilde{D},g \in G}W_{d_{g},g}$ where $d_{g}=0$ except for finitely many $g$. We let $\mu(f):=\displaystyle \sum_{g \in G}\nu(d_{g})w_{g^{-1}}$. We claim that $\mu$ is well-defined. Let $d_{1},d_{2},\cdots,d_{n} \in \widetilde{D}$ and $g_{1},g_{2},\cdots,g_{n} \in G$ (with $g_{i}$'s distinct) be such that $\sum_{i=1}^{n}W_{d_{i},g_{i}}=0$. To show that $\mu$ is well-defined, it is enough to verify that $\sum_{i=1}^{n}\nu(d_{i})w_{g_{i}}^{*}=0$. 

Let $i \in \{1,2,\cdots,n\}$ be given. Let $X \in \Omega \cap \Omega g_{i}^{-1}$ be given. Then $(X,g_{i}) \in \clg$. Now $\sum_{i=1}^{n}W_{d_{i},g_{i}}(X,g_{i})=0$ implies that $d_{i}+I_{X}=0$. Hence $d_{i}+I_{X}=0$ for every $X \in \Omega \cap \Omega g_{i}^{-1}$. Now for $X \in \widetilde{\Omega}$, $1_{\Omega}1_{\Omega g_{i}^{-1}}d_{i}+I_{X}=1_{\Omega}(X)1_{\Omega g_{i}^{-1}}(X)(d_i+I_{X})$. Thus $1_{\Omega}1_{\Omega g_{i}^{-1}}d_{i}+I_{X}=0$ for every $X \in \widetilde{\Omega}$. This implies that $d_{i}1_{\Omega}1_{\Omega g_{i}^{-1}}=0$. Now observe that \[\nu(d_{i})w_{g_{i}^{-1}}=\nu(d_{i})e_{g_{i}^{-1}}w_{g_{i}^{-1}} =\nu(d_{i})\nu(1_{\Omega})\nu(1_{\Omega g_{i}^{-1}})w_{g_{i}^{-1}}=\nu(d_{i}1_{\Omega}1_{\Omega g_{i}^{-1}})w_{g_{i}^{-1}}=0.\]
As a result, the sum $\sum_{i=1}^{n}\nu(d_{i})w_{g_{i}^{-1}}$ vanishes. This proves that $\mu$ is well-defined. 

Next we verify that $\mu$ is multiplicative. We leave it to the reader to convince himself that it is enough to prove that $\mu(W_{d_{1},g_{1}}W_{d_{2},g_{2}})=\mu(W_{d_{1},g_{1}})\mu(W_{d_{2},g_{2}})$ for $d_{1},d_{2} \in \widetilde{D}$ and $g_{1},g_{2} \in G$. 

\textit{Claim:} \label{just proven claim} For $g \in G$ and $d \in \widetilde{D}$, $w_{g}\nu(d)w_{g}^{*}=\nu(\beta_{g}^{-1}(d)1_{\Omega g})$. 

Since the linear span of $\{j_{h}(x): h \in G, x \in A\}$ is dense in $\widetilde{D}$, it is enough to prove the claim when $d=j_{h}(x)$ for some $h \in G$ and $x \in A$. Let $h\in G$ and $x \in A$ be given and let $d:=j_{h}(x)$. Note that $\nu(j_{h}(x))=w_{h}xw_{h}^{*}$. Now calculate as follows to observe that
\begin{align*}
w_{g}w_{h}xw_{h}^{*}w_{g}^{*}&=e_{g}w_{hg}xw_{hg}^{*}e_{g} ~~(\textrm{by Lemma \ref{w-relations}}) \\
                             &=e_{g}\nu(j_{hg}(x)) \\
                             & = e_{g}\nu(\beta_{g}^{-1}(d)) \\
                             &=\nu(1_{\Omega g})\nu(\beta_{g}^{-1}(d)) \\
                             &=\nu(\beta_{g}^{-1}(d)1_{\Omega g})
\end{align*}
This proves the claim. Now let $g_{1},g_{2} \in G$ and $d_{1},d_{2} \in \widetilde{D}$ be given. By Lemma \ref{Relations}, it follows that $W_{d_{1},g_{1}}W_{d_{2},g_{2}}=W_{d,g_{1}g_{2}}$ where $d=1_{\Omega g_{1}^{-1}}d_{1}\beta_{g_{1}}(d_{2})$. Thus by definition, we have $\mu(W_{d_{1},g_{1}}W_{d_{2},g_{2}})=\nu(d)w_{g_{1}g_{2}}^{*}$. Now calculate as follows to observe that 
\begin{align*}
\mu(W_{d_{1},g_{1}})\mu(W_{d_{2},g_{2}})&=\nu(d_{1})w_{g_{1}}^{*}\nu(d_{2})w_{g_{2}}^{*} \\
                                        &=\nu(d_{1})w_{g_{1}}^{*}e_{g_{1}}\nu(d_{2})w_{g_{2}}^{*} \\
                                        &= \nu(d_{1})w_{g_{1}}^{*}\nu(d_{2})w_{g_{1}}w_{g_{1}}^{*}w_{g_{2}}^{*} ~~(\textrm{Since $e_{g}$ commutes with $\mathcal{F}$}) \\
                                        &= \nu(d_{1})w_{g_{1}^{-1}}\nu(d_{2})w_{g_{1}^{-1}}^{*}e_{g_{1}^{-1}}w_{g_{1}g_{2}}^{*} ~~(\textrm{ by Lemma \ref{w-relations}}) \\
                                        &= \nu(d_{1}\beta_{g_1}(d_{2})1_{\Omega g_{1}^{-1}})w_{g_{1}g_{2}}^{*}  \\
                                        &=\nu(d_{1}\beta_{g_1}(d_{2})1_{\Omega g_{1}^{-1}})w_{g_{1}g_{2}}^{*} \\
                                        &= \nu(d)w_{g_{1}g_{2}}^{*}. 
                                        \end{align*} 
Hence $\mu(W_{d_{1},g_{1}}W_{d_{2},g_{2}})=\mu(W_{d_{1},g_{1}})\mu(W_{d_{2},g_{2}})$. This shows that $\mu$ is multiplicative. A similar computation yields that $\mu$ is $*$-preserving. Thus $\mu:\Gamma_{c}(\clg,r^{*}\mathcal{D}) \to A \rtimes P$ is a $*$-homomorphism. Note that by defintion $\mu(W_{d,e})=\nu(d)$ for every $d \in \widetilde{D}$ and the map $\widetilde{D} \ni d \to W_{d,e} \in \Gamma(\Omega, \mathcal{D})$ is a $*$-homomorphism. Thus $\mu$ restricted to $\Gamma(\Omega,\mathcal{D})$ is bounded. Hence the $*$-homomorphism $\mu$ extends to a $*$-homomorphism from the full crossed product $\mathcal{D} \rtimes \clg$ to $A \rtimes P$ which we still denote by $\mu$. 

\begin{thm}
\label{injectivity of lambda}
Let $\lambda:A \rtimes P \to \mathcal{D} \rtimes \clg$ be the $*$-homomorphism constructed in Proposition \ref{homomorphism lambda}. Then $\lambda$ is an isomorphism. 
\end{thm}
\textit{Proof.} Let $\mu: \mathcal{D} \rtimes \clg \to A \rtimes P$ be the $*$-homomorphism such that $\mu(W_{d,g})=\nu(d)w_{g}^{*}$. The existence of such a map is shown in the paragraphs preceeding this theorem. By definition, for $a \in P$, $\mu(\widetilde{v}_{a})=\nu(1_{\Omega})w_{a}=v_{a}$ and for $x \in A$, $\mu(W_{j_{e}(x),x})=\nu(j_{e}(x))w_{e}=x$. Thus $\mu \circ \lambda =id$. This proves that $\lambda$ is $1$-$1$. The surjectivity of $\lambda$ is already proven in Proposition \ref{homomorphism lambda}. Hence $\lambda$ is an isomomorphism. This completes the proof.

 \section{K-group computation for the free semigroup $\mathbb{F}_{n}^{+}$}

As an application of our results, we show that the $K$-theory of $A \rtimes \mathbb{F}_{n}^{+}$ coincides with that of $A$ where $\mathbb{F}_{n}^{+}$ denotes the free semigroup on $n$-generators. For the rest of this paper, let $\mathbb{F}_{n}$ be the free group on $n$ generators and we denote the generators by $a_{1},a_{2},\cdots,a_{n}$. Denote the semigroup generated by $a_{1},a_{2},\cdots a_{n}$ by $\mathbb{F}_{n}^{+}$ i.e. $\mathbb{F}_{n}^{+}$ consists of words in $a_{1},a_{2},\cdots,a_{n}$. We reserve the letters $a_{1},a_{2},\cdots,a_{n}$ to denote the canonical generators of $\mathbb{F}_{n}$ and $\mathbb{F}_{n}^{+}$.
It is due to Nica, ( Example 4, Page 23 of  \cite{Nica92}), that $(\mathbb{F}_{n}^{+},\mathbb{F}_{n})$ is a quasi-lattice ordered pair. We should remark that we consider the right variant of Nica's definition of a quasi-lattice ordered pair. But we can apply Nica's proof by considering  the pair $(Q:=(\mathbb{F}_{n}^{+})^{-1}, \mathbb{F}_{n})$. Note $Q$ is the free semigroup on the generators $a_{1}^{-1},a_{2}^{-1},\cdots a_{n}^{-1}$. Thus it follows that if $gQ \cap Q \neq \emptyset$ for $ g\in \mathbb{F}_{n}$  then there exists $a \in Q$ such that $gQ \cap Q=aQ$. Taking inverses imply that $(\mathbb{F}_{n}^{+},\mathbb{F}_{n})$ is quasi-lattice ordered in our sense. 

Let $\clg$ be the Wiener-Hopf groupoid associated to the pair $(\mathbb{F}_{n}^{+},\mathbb{F}_{n})$. We need the fact that $\clg$ is amenable.
To see this, let $P:=\mathbb{F}_{n}^{+}$ and $Q:=P^{-1}$. For $a \in Q$, let $V_{a}: \ell^{2}(Q) \to \ell^{2}(Q)$ be defined by $V_{a}(\delta_{b})=\delta_{ab}$ for $b \in Q$. Here $\{\delta_{b}: b \in Q\}$ denotes the canonical orthonormal basis of $\ell^{2}(Q)$. The unitary $\ell^{2}(P) \ni \delta_{a} \to \delta_{a^{-1}} \in \ell^{2}(Q)$ induces an isomorphism between the Wiener-Hopf algebra $\mathcal{W}(\mathbb{F}_{n}^{+},\mathbb{F}_{n})$ and the $C^{*}$-algebra generated by $\{V_{a}: a \in P\}$, called the $C^{*}$-algebra associated to the left regular representation of $Q$ and let us denote it by $C_{\ell,r}^{*}(Q)$. By Corollary 8.3 of \cite{Li13}, it follows that the $C_{\ell,r}^{*}(Q)$ is nuclear. Thus the $C^{*}$-algebra $\mathcal{W}(\mathbb{F}_{n}^{+},\mathbb{F}_{n})$ is nuclear. But Theorem \ref{main theorem} implies that $\mathcal{W}(\mathbb{F}_{n}^{+},\mathbb{F}_{n})$ is isomorphic to $C_{red}^{*}(\clg)$. Thus by Theorem 5.6.18 of \cite{Ozawa-Brown}, it follows that $\clg$ is amenable. The following corollary is an immediate consequence of Theorem \ref{main theorem} and Theorem \ref{injectivity of lambda}

\begin{crlre}
\label{universal=red}
Let $A$ be a unital $C^{*}$-algebra and $\alpha:\mathbb{F}_{n}^{+} \to End(A)$ be a  unital left action. The regular representation $\rho:A \rtimes \mathbb{F}_{n}^{+} \to A \rtimes_{red} \mathbb{F}_{n}^{+}$ is an isomorphism.
\end{crlre}
\textit{Proof.} 
Let $(\mathcal{D},\clg,\beta)$ be the Wiener-Hopf groupoid dynamical system associated to $(A,\mathbb{F}_{n}^{+},\mathbb{F}_{n},\alpha)$. Since the Wiener-Hopf groupoid $\clg$ is amenable, it follows from Theorem 1 of \cite{Sims} that the natural map from $\mathcal{D} \rtimes \clg$ to the reduced crossed product $\mathcal{D} \rtimes_{red} \clg$ is an isomorphism. The proof is complete by applying Theorem \ref{main theorem} and Theorem \ref{injectivity of lambda}. \hfill $\Box$

Let $A$ be a $C^*$-algebra and $\alpha:\mathbb{F}_{n}^{+} \to End(A)$ be a left action. For $i \in \{1,2,\cdots,n\}$, let $\alpha_{i}:=\alpha_{a_{i}}$. Note that $\alpha$ is completely determined by the $n$ endomorphisms $\alpha_{1},\alpha_{2},\cdots ,\alpha_{n}$. Conversely let  $\alpha_{1},\alpha_{2},\cdots, \alpha_{n}$ be endomorphisms of $A$. Then there exists a unique action $\alpha:\mathbb{F}_{n}^{+} \to End(A)$ such that $\alpha_{a_{i}}=\alpha_{i}$. To see this, let $w:=a_{i_1}a_{i_2}\cdots a_{i_k}$ be a word in $\{a_{1},a_{2},\cdots,a_{n}\}$. Since any word in $\{a_{1},a_{2},\cdots,a_{n}\}$ is a reduced word [See \cite{Serre}] , it follows that the expression $a_{i_{1}}a_{i_{2}}\cdots a_{i_{k}}$ representing $w$ is unique. Now set $\alpha_{w}=\alpha_{i_1}\alpha_{i_2}\cdots \alpha_{i_k}$.  Then $\alpha:\mathbb{F}_{n}^{+} \to End(A)$ is the required action. 

Let $\clh$ be a Hilbert space and let $V:\mathbb{F}_{n}^{+} \to B(\clh)$ be an antihomomorphism such that for each $a \in \mathbb{F}_{n}^{+}$, $V_{a}$ is an isometry. For $i \in \{1,2,\cdots,n\}$, let $V_{i}:=V_{a_{i}}$. Then $V_{1},V_{2},\cdots,V_{n}$ are isometries. Conversely, let $V_{1},V_{2},\cdots, V_{n}$ be isometries on $\clh$. Then there exists a unique antihomorphism $V:\mathbb{F}_{n}^{+} \to B(\clh)$ such that $V_{a_{i}}=V_{i}$ for $i=1,2,\cdots,n$. We leave this verification to the reader. 

\begin{lmma}
\label{Nica-covariance}
Let $A$ be a $C^{*}$-algebra and $\alpha:\mathbb{F}_{n}^{+} \to End(A)$ be a left action. Let $\pi:A \to B(\clh)$ be a $*$-homomorphism and $V:\mathbb{F}_{n}^{+} \to B(\clh)$ be an antihomomorphism of isometries. For $i=1,2,\cdots,n$, let $\alpha_{i}:=\alpha_{a_{i}}$ and $V_{i}:=V_{a_{i}}$. Then the following are equivalent.
\begin{enumerate}
\item[(1)] The pair $(\pi,V)$ is Nica covariant.
\item[(2)] For $i \in \{1,2,\cdots,n\}$ and $x \in A$, $\pi(x)V_i=V_{i}\alpha_{i}(x)$.  For $i,j \in \{1,2,\cdots,n\}$ with $i\neq j$, $V_{i}^{*}V_{j}=0$. 
\end{enumerate}
\end{lmma}
\textit{Proof.} Note that for $i \neq j$, $\mathbb{F}_{n}^{+}a_{i} \cap \mathbb{F}_{n}^{+}a_{j} =\emptyset$. Thus clearly $(1)$ implies $(2)$. Now assume that $(2)$ holds. For $w \in \mathbb{F}_{n}^{+}$, let $E_{w}:=V_{w}V_{w}^{*}$. We leave it to the reader to verify that $\pi(x)V_{w}=V_{w}\alpha_{w}(x)$ for $x \in A$ and $w \in \mathbb{F}_{n}^{+}$.

Observe that for $w_{1},w_{2} \in \mathbb{F}_{n}^{+}$, 
\begin{equation*}
\mathbb{F}_{n}^{+}w_1 \cap \mathbb{F}_{n}^{+}w_2 = \begin{cases}
              \mathbb{F}_{n}^{+}w_2 & \mbox{if }  w_1 \leq w_2, \cr
              &\cr
              \mathbb{F}_{n}^{+}w_1 &  \mbox{if }  w_2 \leq w_1, \cr
              &\cr
              \emptyset & \mbox{else.}  \cr      
\end{cases}
\end{equation*}
Let $w_{1},w_{2} \in \mathbb{F}_{n}^{+}$ be given. Suppose $w_{1} \leq w_{2}$. Let $w \in \mathbb{F}_{n}^{+}$ be such that $ww_{1}=w_{2}$. Note that $E_{w_{2}}=V_{ww_{1}}V_{ww_{1}}^{*}=V_{w_{1}}(V_{w}V_{w}^{*})V_{w_1} \leq V_{w_1}V_{w_1}^{*}=E_{w_{1}}$. Thus $E_{w_{1}}E_{w_{2}}=E_{w_{2}}$ if $w_{1} \leq w_{2}$. Thus  if $w_{1},w_{2},w \in \mathbb{F}_{n}^{+}$ are such that $\mathbb{F}_{n}^{+}w_{1} \cap \mathbb{F}_{n}^{+}w_{2}=\mathbb{F}_{n}^{+}w$ then $E_{w_1}E_{w_{2}}=E_{w}$. Now let $w_{1},w_{2} \in \mathbb{F}_{n}^{+}$ be such that $\mathbb{F}_{n}^{+}w_{1} \cap \mathbb{F}_{n}^{+}w_{2} = \emptyset$. Write $w_1=a_{i_1}a_{i_2}\cdots a_{i_n}$ and $w_2=a_{j_1}a_{j_2}\cdots a_{j_m}$ with $i_{k}, j_{\ell} \in \{1,2,\cdots,n\}$. We claim that $E_{w_1}E_{w_2}=0$. Without loss of generality, we can assume that $n \leq m$. Let $k \in \{0,1,\cdots,n-1\}$ be the least integer for which $i_{n-k} \neq j_{m-k}$. Such an integer exists, otherwise $w_{1} \leq w_{2}$ which contradicts the assumption that the intersection $\mathbb{F}_{n}^{+}w_{1} \cap \mathbb{F}_{n}^{+}w_{2}$ is empty. 

Now use the fact that $V_{i}$ is an isometry to observe that 
\begin{align*}
V_{w_{1}}^{*}V_{w_{2}}&=V_{i_1}^{*}V_{i_{2}}^{*}\cdots V_{i_{n-k}}^{*}V_{j_{m-k}}V_{j_{m-k-1}}\cdots V_{j_{1}}\\
                      &=0 \textrm{~~\big( Since $V_{i}^{*}V_{j}=0$ if $i \neq j$ and $i_{n-k} \neq j_{m-k}$ \big)}. 
\end{align*}
Thus $E_{w_{1}}E_{w_{2}}=0$ if $\mathbb{F}_{n}^{+}w_{1} \cap \mathbb{F}_{n}^{+}w_{2} = \emptyset$. This proves that $(\pi,V)$ is Nica covariant. This completes the proof. \hfill $\Box$

Let $A$ be a unital $C^{*}$-algebra and $\alpha:\mathbb{F}_{n}^{+} \to End(A)$ be a unital left action.  We show that the natural inclusion  $A \ni x \to x \in A \rtimes  \mathbb{F}_{n}^{+}$ induces an isomorphism at the $K$-theory level. Let $(\pi,V)$ be the standard Nica covariant pair for $(A,\mathbb{F}_{n}^{+},\alpha)$. 

 For $i=1,2,\cdots,n$, let $\alpha_{i}:=\alpha_{a_{i}}$.  Consider the Hilbert $A$-modules $\mathcal{E}^{(0)}=A \otimes \ell^{2}(\mathbb{F}_{n}^{+})$ and $\mathcal{E}^{(1)}=A \otimes \ell^{2}(\mathbb{F}_{n}^{+}\backslash \{1\})$. Here $1$ denote the identity element of $\mathbb{F}_{n}^{+}$.  Note that for $x \in A$ and $a \in \mathbb{F}_{n}^{+}$, $\pi(x)$ and $V_{a}$ leaves $\mathcal{E}^{(1)}$ invariant. For $x \in A$ and $a \in \mathbb{F}_{n}^{+}$, denote  the restriction of $\pi(x)$ and $V_{a}$ on $\ell^{2}(\mathbb{F}_{n}^{+}\backslash \{1\})$ by $\widetilde{\pi}(x)$ and $\widetilde{V}_{a}$. Note that $(\widetilde{\pi},\widetilde{V})$ is Nica covariant. The only things that needs verification is that  for $i \neq j$, the range of $\widetilde{V}_{a_{i}}$ is orthogonal to that of $V_{a_{j}}$. This follows from the fact that the range of $V_{a_{i}}$ is $A \otimes \ell^{2}((\mathbb{F}_{n}^{+}\backslash \{1\})a_{i})$.   Let us denote the map $\pi \rtimes V: A \rtimes \mathbb{F}_{n}^{+} \to \mathcal{L}_{A}(\mathcal{E}^{(0)})$ by $\lambda^{(0)}$ and the map $\widetilde{\pi} \rtimes \widetilde{V}: A \rtimes \mathbb{F}_{n}^{+} \to \mathcal{L}_{A}(\mathcal{E}^{(1)})$ by $\lambda^{(1)}$ respectively. Denote the orthogonal projection from $\ell^{2}(\mathbb{F}_{n}^{+})$ onto $\ell^{2}(\mathbb{F}_{n}^{+}\backslash \{1\})$ by $Q$ and  Let $P:\mathcal{E}^{(0)} \to \mathcal{E}^{(1)}$ be defined by $P:=1 \otimes Q$. 

\textbf{Notations:} Let $\{\delta_{x}:x \in \mathbb{F}_{n}^{+}\}$ be the standard orthonormal basis for $\ell^{2}(\mathbb{F}_{n}^{+})$ and let $\{e_{x,y}:x, y \in \mathbb{F}_{n}^{+}\}$ be the standard 'matrix units' with respect to the orthonormal basis $\{\delta_{x}:x \in \mathbb{F}_{n}^{+}\}$. For $i \in \{1,2,\cdots,n\}$, let $v_{i}:\ell^{2}(\mathbb{F}_{n}^{+}) \to \ell^{2}(\mathbb{F}_{n}^{+})$ be given by $v_{i}(\delta_{a})=\delta_{aa_{i}}$. Here $\{\delta_{a}:a \in \mathbb{F}_{n}^{+}\}$ stands for the standard orthonormal basis. Let $p$ denote the orthogonal projection from $\ell^{2}(\mathbb{F}_{n}^{+})$ onto the one dimensional space subspace $\ell^{2}(\{1\})$ and set $q:=1-p$.

\begin{lmma}
\label{KK-element}
The triple $\Big(\mathcal{E}:=\mathcal{E}^{(0)} \oplus \mathcal{E}^{(1)}, \lambda:= \lambda^{(0)} \oplus \lambda^{(1)}, F:= \begin{bmatrix}
                                                                                                           0 & P^{*} \\
                                                                                                           P & 0
                                                                                                             \end{bmatrix} \Big) \in KK(A \rtimes \mathbb{F}_{n}^{+},A)$.
\end{lmma}
\textit{Proof.} Note that $PP^{*}=1$ and $P^{*}P=1-1\otimes p$. Thus $P^{*}P-1$ and $PP^{*}-1$ are compact. Note that we have assumed that $A$ is unital. 

Note that for $x \in A$ and $a,b \in \mathbb{F}_{n}^{+}$, 
\begin{align*}
P\pi(x)-\widetilde{\pi}(x)P&=0, \\
PV_{a}-\widetilde{V}_{a} P&= 1 \otimes pe_{a,1}, \textrm{~and} \\
PV_{b}^{*}-\widetilde{V}_{b}^{*}P&=0.
\end{align*}

From the above equations and the fact that the linear span of $\{v_{a}xv_{b}^{*}:a ,b \in \mathbb{F}_{n}^{+}, x \in A\}$ is dense in $A \rtimes \mathbb{F}_{n}^{+}$, it follows that $P\lambda^{(0)}(T)-\lambda^{(1)}(T)P$ is compact for every $T \in A \rtimes \mathbb{F}_{n}^{+}$.  The proof is now complete. \hfill $\Box$ 

Let us denote the Kasparov triple defined in Lemma \ref{KK-element} by $[d]$. Let us denote the inclusion $A \ni x \to x \in A \rtimes \mathbb{F}_{n}^{+}$ by $j$ and denote the corresponding $KK$-element representing $j$ in $KK(A, A \rtimes \mathbb{F}_{n}^{+})$ by $[j]$. We use the notation $\sharp$ to denote the Kasparov product. We adapt the proof of Theorem 2.3 of \cite{KS97}.

\begin{thm}
\label{K-theory}
The $KK$-elements $[d]$ and $[j]$ are inverses of each other with respect to the Kasparov product.
\end{thm}
\textit{Proof.} Note the decomposition $\mathcal{E}^{(0)}=A \oplus \mathcal{E}^{(1)}$ as Hilbert $A$-modules. With respect to this decomposition, observe that, the pull-back of $[d]$ by the homomorphism $j$ i.e. the Kasparov product $[j] \sharp [d]$ is isomorphic to the direct sum of  $(A,m,0)$ and $\mathcal{D}:=\Big(\mathcal{E}^{(1)} \oplus \mathcal{E}^{(1)}, \widetilde{\pi} \oplus \widetilde{\pi}, \begin{bmatrix}
                                                                                                             0 & 1 \\
                                                                                                             1 & 0
                                                                                                             \end{bmatrix} \Big)$ where $m:A \to A=M(A)$ is the usual multiplication representation. Observe that $\mathcal{D}$ is degenerate. Thus the product $[j] \sharp [d] =[1_{A}]$. 
                                                                                                              
For this proof, let us denote $A \rtimes \mathbb{F}_{n}^{+}$ by $\mathcal{T}$.  Note that the Kasparov product $[d] \sharp [j]$, the push-forward of $[d]$ by the homomorphism $j$, is given by $(\mathcal{E} \otimes_{A} \mathcal{T},\sigma:= \lambda \otimes 1,G:= F \otimes 1)$. The map \[(\ell^{2}(\Gamma) \otimes A)\otimes_{A} \mathcal{T}  \ni (\delta_{a} \otimes x) \otimes y\to \delta_{a} \otimes \pi(x)y \in  \ell^{2}(\Gamma) \otimes \mathcal{T}\] is  unitary (The surjectivity follows from the assumption that the action $\alpha$ is unital and hence $\pi$ is unital). Here $\Gamma$ stands for either $\mathbb{F}_{n}^{+}$ or $\mathbb{F}_{n}^{+} \backslash \{1\}$.  

We identify this way the Hilbert $\mathcal{T}$-modules, $\mathcal{E}^{(0)} \otimes_{A} \mathcal{T}$ with $\mathcal{F}^{(0)}:=\ell^{2}(\mathbb{F}_{n}^{+}) \otimes \mathcal{T}$ and $\mathcal{E}^{(1)}\otimes_{A} \mathcal{T}$ with $\mathcal{F}^{(1)}:=\ell^{2}(\mathbb{F}_{n}^{+}\backslash \{1\}) \otimes \mathcal{T}$. With this identification, the representation $\sigma=\sigma^{(0)} \oplus \sigma^{(1)}$ of $\mathcal{T}$ is given by the formulas: for $x \in A$ and $a \in \mathbb{F}_{n}^{+}$,
\begin{align*}
\sigma^{(i)}(x)(\delta_{c} \otimes y)&:=\delta_{c} \otimes \alpha_{c}(x)y \\
\sigma^{(i)}(v_{a})(\delta_{c} \otimes y)&:=\delta_{ca} \otimes y 
\end{align*}
Let $Q$ be the orthogonal projection of $\ell^{(2)}(\mathbb{F}_{n}^{+})$ onto $\ell^{(2}(\mathbb{F}_{n}^{+} \backslash \{1\})$ and let $R:\mathcal{F}^{(0)} \to \mathcal{F}^{(1)}$ be given by $R:=Q \otimes 1$. The operator $G=F \otimes 1$ is then given by $\begin{bmatrix}
0 & R^{*} \\
R & 0
\end{bmatrix}.$

For $i=1,2,\cdots,n$ and $t \in [0,\frac{\pi}{2}]$, let $w_{i}^{(t)} \in \mathcal{L}_{\mathcal{T}}(\ell^{2}(\mathbb{F}_{n}^{+}) \otimes \mathcal{T})$ be defined by 
\[
w_{i}^{(t)}:=\cos(t)(v_ip \otimes 1) + \sin(t)(p \otimes V_{i})+v_{i}(1-p) \otimes 1
\]
where $V_{i}=V_{a_{i}}$ and $p$ denotes the orthognal projection of $\ell^{2}(\mathbb{F}_{n}^{+})$ onto the one-dimensional subspace spanned by $\{\delta_{1}\}$. Observe that $\{w_{i}^{t}:i=1,2,\cdots,n\}$ is a collection of isometries with orthogonal range projections. For $t \in [0,\frac{\pi}{2}]$, let $w^{t}:\mathbb{F}_{n}^{+} \to \mathcal{L}_{\mathcal{T}}(\mathcal{F}^{(0)})$ be the unique antihomomorphism of isometries such that $w^{t}_{a_{i}}=w^{t}_{i}$. 

For $t \in [0,\frac{\pi}{2}]$, let $\pi^{t}:A \to \mathcal{L}_{\mathcal{T}}(\mathcal{F}^{(0)})$ be defined by $\pi^{t}(x)=\sigma^{(0)}(x)$ for $x \in A$. Note that for $t \in [0,\frac{\pi}{2}]$, $x \in A$ and $i\in \{1,2,\cdots,n\}$, $\pi^{t}(x)w_{i}^{t}=w_{i}^{t}\pi^{t}(\alpha_{i}(x))$. By Lemma \ref{Nica-covariance}, it follows that the pair $(\pi^{t},w^{t})$ is Nica-covariant. Thus there exists a homomorphism $\sigma^{(t)}:\mathcal{T} \to \mathcal{L}_{\mathcal{T}}(\mathcal{F}^{(0)})$ such that $\tau^{(t)}(x)=\pi^{t}(x)$ and $\tau^{(t)}(v_{a})=w_{a}^{t}$. Observe that $\tau^{(0)}=\sigma^{(0)}$. Also note that for $t \in [0,\frac{\pi}{2}]$, $\tau^{(t)}(v_{a_{i}})-\sigma^{(0)}(v_{a_{i}})$ is compact for every $i=1,2,\cdots,n$ and $\tau^{(t)}(x)=\sigma^{(0)}(x)$ for $x \in A$. Since $\mathcal{T}$ is generated by $A$ and the isometries $\{v_{a_{i}}:i=1,2,\cdots,n\}$, it follows that for every $s \in \mathcal{T}$, $\tau^{(t)}(s)-\sigma^{(0)}(s)$ is compact. 

Hence in $KK(\mathcal{T},\mathcal{T})$, $[d] \sharp [j]=(\mathcal{F}^{(0)}\oplus \mathcal{F}^{(1)}, \tau^{(t)}\oplus \sigma^{(1)}, G)$ for every $t \in [0,\frac{\pi}{2}]$. Observe that the decomposition $\mathcal{F}^{(0)}=\mathcal{T} \oplus \mathcal{F}^{(1)}$ is left invariant by
$\tau^{\frac{\pi}{2}}$. With respect to this decomposition $\tau^{\frac{\pi}{2}}=m \oplus \sigma^{(1)}$ where $m:\mathcal{T} \to \mathcal{T}$ is the identity map. Thus $[d] \sharp [j] =(\mathcal{T}, m, 0) \oplus \mathcal{D}$ where $\mathcal{D}=(\mathcal{F}^{(1)} \oplus \mathcal{F}^{(1)}, \sigma^{(1)}\oplus \sigma^{(1)}, \begin{bmatrix}
0 & 1 \\
1 & 0
\end{bmatrix})$. Note that $\mathcal{D}$ is degenerate. As a consequence, it follows that $[d] \sharp [j]=[1_{\mathcal{T}}]$ in $KK(\mathcal{T},\mathcal{T})$. This completes the proof. \hfill $\Box$

\begin{rmrk}
Let $A$ be a $C^{*}$-algebra and $\alpha:\mathbb{F}_{n}^{+} \to End(A)$ be a left action. We do not assume that $A$ is unital or $\alpha$ is unital. Theorem \ref{K-theory} and a six term exact sequence argument applied to the short exact sequence of Lemma \ref{unitisation} implies that $K_{*}(A) \cong K_{*}(A \rtimes_{red} \mathbb{F}_{n}^{+})$. 
\end{rmrk}

 \bibliography{references}
 
\bibliographystyle{amsalpha}

\noindent
{\sc S. Sundar}
(\texttt{sundarsobers@gmail.com})\\
         {\footnotesize  Chennai Mathematical Institute, H1 Sipcot IT Park, \\
Siruseri, Padur, 603103, Tamilnadu, INDIA.}

\end{document}